\documentclass[12pt]{amsart}
\usepackage{amsmath,amsxtra,amssymb,amsfonts}

\newcommand{\8}{\infty}

\newcommand{\kla}{\left ( }
\newcommand{\mer}{\right ) }

\newcommand{\for}{\begin{eqnarray*}}
\newcommand{\mel}{\end{eqnarray*}}

\newcommand{\mitt}{\left | { \atop } \right.}
\newcommand{\kl}{\pl \le \pl}

\newcommand{\lel}{\pl = \pl}

\newcommand{\nz}{{\mathbb N}}
\newcommand{\nen}{n \in \nz}

\newcommand{\rz}{{\mathbb R}}

\newcommand{\cz}{{\mathbb C}}

\newcommand{\ten}{\otimes}

\newcommand{\p}{\hspace{.05cm}}
\newcommand{\pl}{\hspace{.1cm}}
\newcommand{\pll}{\hspace{.3cm}}

\newcommand{\hz}{\vspace{0.5cm}}

\renewcommand{\qed}{\hspace*{\fill} {\vrule height7pt width7pt 
depth0pt} \hz\pagebreak[1]}

\newcommand{\Om}{\Omega}
\newcommand{\om}{\omega}
\renewcommand{\a}{\alpha}
\newcommand{\al}{\alpha}
\newcommand{\si}{\sigma}
\newcommand{\Si}{\Sigma}

\newcommand{\la}{\lambda}
\newcommand{\eps}{\varepsilon}

\renewcommand{\L}{{\mathcal L}}

\newcommand{\E}{{\mathcal E}}
\newcommand{\V}{{\mathcal L}}
\newcommand{\A}{{\mathcal A}}
\newcommand{\B}{{\mathcal B}}

\newcommand{\M}{{\mathcal M}}
\newcommand{\K}{{\mathcal K}}

\newcommand{\N}{{\mathcal N}}

\newcommand{\hhz}{\vspace{0.3cm}}

\newcommand{\noo}{\left \|}
\newcommand{\rrm}{\right \|}

\newcommand{\bet}{\left |}
\newcommand{\rag}{\right |}

\newcommand{\intt}{\int\limits}
\newcommand{\summ}{\sum\limits}

\newcommand{\ez}{{\rm I\!E}}

\newtheorem{lemma}{Lemma}[section]
\newtheorem{prop}[lemma]{Proposition}
\newtheorem{theorem}[lemma]{Theorem}
\newtheorem{cor}[lemma]{Corollary}

\newtheorem{rem}[lemma]{Remark}

\newcommand{\re}{\begin{rem}\rm}
  \newcommand{\mar}{\end{rem}}

\newtheorem{defi}[lemma]{Definition}
\renewcommand{\baselinestretch}{1.2}
\begin{document}

\author[Marius Junge]{Marius Junge$^{\dag}$}
\address{Department of Mathematics\\
University of Illinois at  Urbana-Champaign, Urbana,
IL 61801} \email{junge\@math.uiuc.edu}

\thanks{$^{\dag}$Junge is  partially supported by the National Science
Foundation}

\title[Non-commutative  Doob inequality]
{\bf \large Doob's inequality for non-commutative martingales}




\hz

\begin{abstract}
Let $1\le p<\8$ and $(x_n)_{\nen}$ be a sequence of positive
elements in a non-commutative $L_p$ space and $(E_n)_{\nen}$ be an
increasing sequence of conditional expectations, then
\[ \noo \summ_n E_n(x_n) \rrm_p \kl c_p \noo \summ_n x_n \rrm_p \pl. \]
This inequality is due to Burkholder, Davis and Gundy in the
commutative case. By duality, we obtain a version of Doob's
maximal inequality for  $1<p\le \8$.
\end{abstract}

\maketitle

{\bf Introduction:}

Inspired by  quantum mechanics and probability, non-commutative
probability has become an independent field of mathematical
research. We refer to P.A. Meyer's exposition \cite{M}, the
successive conferences on quantum probability \cite{AvW}, the
lecture notes by Jajte \cite{Ja1,Ja2} on almost sure and uniform
convergence and finally   the work of Voiculescu,  Dykema, Nica
\cite{VDN} and of Biane, Speicher \cite{BS} concerning the
recent  progress in free probability and free Brownian motion.
Doob's inequality is a classical tool in probability and
analysis. Transferring  classical inequalities into the
non-commutative setting theory often requires an additional
insight. Pisier, Xu \cite{PX,Pcomb} use functional analytic  and
combinatorial methods to establish the non-commutative versions
of the Burkholder-Gundy square function inequality. The absence
of stopping time arguments, at least until the time of this
writing, imposes one of the main difficulties in this recent
branch of martingale theory.

\hz

The formulation of Doob's inequality for non-commutative
martingales faces the following  problem. For an increasing
sequence of conditional expectations $(E_n)_{\nen}$ and a
positive operator $x$ in $L_p$, there is no reason why $\sup_n
E_n(x)$ or $\sup_n |E_n(x)|$ should be an element in $L_p$ or
represent a (possibly unbounded) operator at all. Using Pisier's
non-commutative vector-valued $L_p$-space $L_p(N;\ell_\8)$ we can
overcome this problem and  at least guess the right formulation
of Doob's inequality. However, Pisier's definition is  restricted
to von Neumann algebras with a $\si$-weakly dense net of finite
dimensional subalgebras, so-called hyperfinite von Neumann
algebras. But  maximal inequalities are also interesting for free
stochastic processes where the underlying von Neumann algebra is
genuinely  not hyperfinite. All these obstacles disappear for the
so-called dual version of Doob's inequality: For every sequence
$(x_n)_{\nen}$ of positive operators
 \for
  \noo \summ_n E_n(x_n)\rrm_p &\le&  c_p \pl \noo \summ_n x_n\rrm_p  \pl . \\[-1.66cm]
  \mel  \hspace*{\fill} $(DD_p)$

\vspace{0.75cm}

In the commutative case this inequality is due to Burkholder,
Davis and Gundy \cite{B1} (even in the more general setting of
Orlicz norms). Since it is crucial  to understand  our approach
to Doob's inequality, let us indicate the duality argument
relating $(DD_p)$ and Doob's inequality in the commutative case.
Indeed, $(DD_p)$ implies that $T(x_n)=\sum_n E_n(x_n)$ defines a
continuous linear map between $L_p(\ell_1)$ and $L_p$. The norm
of $\noo T^*\rrm$ yields the best constant in Doob's inequality
for the conjugate index $p'=\frac{p}{p-1}$.
  \[ \noo \sup_n |E_n(x)|\rrm_{p'} \lel \noo T^*(x)\rrm_{L_{p'}(\ell_\8)}
 \kl \noo T\rrm \pl   \noo x\rrm_{p'}  \lel c_p \pl \noo x\rrm_{p'}  \pl .\]
Personally, I learned this argument after reading
Dilworth's paper \cite{Dil}. But I am sure it is
known to  experts in the field, see Garcia's
\cite{Gar} for the general theory (and \cite{MS}
for the explicit equivalence). $(DD_p)$ admits an
entirely elementary proof in the commutative case
(see again \cite{MS}). This elementary proof still
works in the non-commutative case for $p=2$, see
Lemma \ref{p2}. It is the starting point of our
investigation. We recommend the reader (not
familiar with modular theory) to start in section 3
where interpolation is used to extend $(DD_p)$ to
$1\le p\le 2$ and suitable norms are introduced to
make the above duality argument work in the
non-commutative case. In section 4, we establish
the dual version $(DD_p)$ in the range $2\le p<\8$
using duality arguments which rely on Pisier/Xu's
version of Stein's inequality in combination with
techniques from Hilbert $C^*$-modules. By duality,
we obtain the non-commutative Doob inequality in
the more delicate range $1<p\le 2$. The heart of our
arguments rely on the (apparently new) connection
between Hilbert $C^*$-modules and non-commutative
$L_p$ spaces presented in section 2. These duality
techniques are necessary because $p>1$ and  $0\le
a\le b$ implies $0\le a^p\le b^p$ only for
\underline{commuting} operators. This is very often
used in the 'elementary' approach to commutative
martingales inequalities as in Garsia's book
\cite{Gar}.

\hz

Let us formulate our main results for finite von Neumann algebras.
If $\tau:N\to \cz$ is a normal, tracial state, i.e.
$\tau(xy)=\tau(yx)$, then the space $L_p(N,\tau)$ is defined by
the completion of $N$ with respect to the norm
 \[ \noo x\rrm_p \lel \tau( (x^*x)^{\frac{p}{2}})^{\frac{1}{p}} \pl .\]
We refer to the first section for more precise definitions and
references. Given a subalgebra $M\subset N$, the embedding
$\iota:L_1(M,\tau)\subset L_1(N,\tau)$ is isometric because
$|x|=\sqrt{x^*x}\in M$ for all $x\in M$. The dual map
$E=\iota^*:N\to M$ yields a conditional expectation satisfying
 \[ E(axb) \lel aE(x)b \]
for all $a,b\in M$ and $x\in N$, see \cite[Theorem 3.4.]{TAK}.
Since $E$ is trace preserving,  $E$ extends to a contraction
$E:L_p(N,\tau)\to L_p(N,\tau)$ with range $L_p(M,\tau)$. In the
following, we consider an increasing sequence
$(N_n)_{\nen}\subset N$ of von Neumann subalgebras with
conditional expectations $(E_n)_{\nen}$. We recall that an
element $x$ is positive if it is of the form $x=y^*y$. \hhz

\begin{theorem}\label{0.1.}\label{ep}  Let $1\le p<\8$, then there exists a  constant
$c_p$ depending only on $p$ such that for every sequence of
positive elements $(x_n)_{\nen} \subset L_p(N,\tau)$
 \for
  \noo \summ_n E_n(x_n)\rrm_p &\le&  c_p \pl \noo \summ_n x_n\rrm_p  \pl . \\[-1.66cm]
  \mel  \hspace*{\fill} $(DD_p)$\hz
\end{theorem}
\hhz

Note the close relation to the non-commutative Stein inequality,
see \cite[Theorem 2.3.]{PX},
 \[ \noo \summ_n  E_n(x_n)^*E_n(x_n)\rrm_{p}  \kl  \gamma_{2p}^2  \pl \noo \summ_n x_n^*x_n\rrm_{p}   \pl .\]
Using Kadison's inequality $E_n(x_n)^*E_n(x_n)\le
E_n(x_n^*x_n)$, it turns out that $(DD_p)$ is
stronger than Stein's inequality. However, Stein's
inequality combined with the theory of Hilbert
$C^*$-modules yields  one of the fundamental
inequalities in the proof of Theorem \ref{0.1.}.
Using a Hahn-Banach separation argument \`a la
Grothendieck-Pietsch, we deduce Doob's maximal
inequality.\hhz

\begin{theorem}\label{0.2.}{\rm [{\bf Doob's maximal inequality}\rm]} Let
$1<p\le \8$ and $x\in L_p(N,\tau)$, then there
exist $a,b\in L_{2p}(N,\tau)$ and a sequence of
contractions $(y_n)\subset N$ such that
 \[ E_n(x) \lel ay_nb  \quad \mbox{and}\quad \noo a\rrm_{2p} \noo b\rrm_{2p} \kl c_{p'} \noo x\rrm_{p} \pl \]
Here $c_{p'}$ is the constant in Theorem \ref{0.1.}
for the conjugate index $p'=\frac{p}{p-1}$. In
particular, for every positive $x\in L_p(N,\tau)$,
there exists a positive $b\in L_p(N,\tau)$ such
that
 \[ E_n(x)\kl b\]
for all $\nen$.
\end{theorem}\hhz

In the case of hyperfinite von Neumann algebras
this is equivalent to the corresponding
vector-valued inequality, see \cite{Pvp}, and
therefore justifies the name 'maximal inequality'.
In the semi-commutative case where
$N_n=L_\8(\Om,\Si_n,\mu)\ten M$ and $\Si_n$ are
increasing $\si$-algebras this inequality is
stronger than the vector-valued Doob inequality,
see Remark \ref{scc}. In particular, this applies
for random matrices. The inequality can be extended
to a continuous index set under suitable density
assumptions, for examples for Clifford martingales
or free stochastical processes. Since these
modifications are rather obvious, we omit  the
details.

\hz

Clearly, maximal inequalities immediately imply almost sure
convergence. Therefore it is not surprising that Theorem
\ref{0.2.}  implies almost uniform convergence convergence of the
martingale truncations  $(E_n(x))$ for $p\ge 2$ and bilateral
convergence for $1<p\le 2$ in case  of a tracial state. We refer
to \cite{Ja1,Ja2} for  the definition of these notions and more
details.  In the tracial case the bilateral convergence of the
martingale truncations is known by a result of Cuculescu
\cite{Cu} even for martingales in $L_1(N,\tau)$. Therefore
Theorem \ref{0.2.} provides a alternative approach to these
results but only for $p>1$. However, the maximal inequality
discussed in \cite{Ja1} cannot easily be interpolated to obtain
Theorem \ref{0.2.} as in the real case. In a subsequent paper
\cite{DJ2}, we will apply the maximal inequality of Theorem
\ref{0.2.} in Haagerup $L_p$ spaces in order to  obtain
(bilateral) almost sure  convergence for all states thus
underlining the strength of these maximal inequalities.

\hz

Preliminary results and notation are contained in section 1.
Section 5 contains immediate applications to submartingales and
conditional expectations associated to actions of groups.

\hz

I am indebted to Q. Xu for many  discussion and support. I want
to thank A. Defant for the discussion about almost everywhere
convergence. The knowledge of similar results for almost  sure
convergence of unconditional sequences, see \cite{DJ}, have been
very encouraging. I would like to thank Stanislaw Goldstein for 
initiating a correction in the proofs of Lemma 2.3 and Lemma 3.2.

\section{\bf Notation and preliminary results}

As a shortcut, we use $(x_n)$, $(x_{nk})$ instead of
$(x_n)_{\nen}$, $(x_{nk})_{n,k\in \nz}$ for sequences indexed by
the natural numbers $\nz$ or its cartesian product $\nz^2$. We use
standard notation in operator algebras, as in \cite{TAK,Kad}. In
particular, $B(H)$ denotes the algebra of bounded operators on a
Hilbert space $H$ and $\K(H)$, $\K$ denote the subalgebra of
compact operators on $H$, $\ell_2$, respectively. The letters $N$,
$M$ will be used for von Neumann algebras, i.e. subalgebras of
some $B(H)$ which are closed with respect to the $\si$-weak
operator topology. We refer to \cite{DX,TAK} for the different
locally convex topologies relevant to operator algebras. For
$\nen$ we denote by $M_n(N)$ the von Neumann algebra of $n\times
n$ matrices with values in $N$. We will briefly use $M_n$ for
$M_n(\cz)$. Given $C^*$-algebras $\A$, $\B$, we denote by $\A \ten
\B$ the minimal tensor product. For von Neumann algebras $N\subset
B(H_1)$, $M\subset B(H_2)$, we use $N\bar{\ten}M$ for the closure
of $N \ten M\subset B(H_1\ten H_2)$ in the $\si$-weak operator
topology. Let us recall that a von Neumann algebra is semifinite
if there exists a normal, semifinite faithful trace. A trace is a
positive homogeneous, additive function on $N_+=\{x^*x \p | \p
x\in N\}$, the cone of positive elements of $N$, such that for
all increasing nets $(x_i)_i$ with supremum in $N$ and for all
$x\in N$
\begin{enumerate}
 \item[n)] $\tau(\sup_i x_i)=\sup_i \tau(x_i)$.
 \item[s)] For every $0<x$ there exists $0<y<x$ such that $0<\tau(y)<\8$.
 \item[f)] $\tau(x)=0$ implies $x=0$.
 \item[t)] For all unitaries $u\in N$: $\tau(uxu^*)=\tau(x)$.
\end{enumerate}
A positive homogeneous, additive function $w:N_+\to [0,\8]$
satisfying n), s), f), but not the last property t), is called a
n.s.f. (normal, semifinite, faithful) weight.

\hz

It will be worthwhile to clarify the different notions of
non-commutative $L_p$-spaces. If $\tau$ is a  trace then
 \[ m(\tau)\lel \left\{\summ_{i=1}^n  y_ix_i \mitt \nen, \summ_{i=1}^n [\tau(y_i^*y_i) + \tau(x_i^*x_i)] <\8 \right\}   \]
is the definition ideal on which there exists a unique linear
extension $\tau:m(\tau) \to \cz$  satisfying $\tau(xy)=\tau(yx)$.
The $L_p$-(quasi)-norm is defined for $x\in m(\tau)$ by
 \[ \noo x\rrm_p \lel \tau( (x^*x)^{\frac p2})^{\frac 1p} \pl .\]
Then $L_p(N,\tau)$  is the completion of $m(\tau)$
with respect to the $L_p$-norm. (For $p<1$  smaller
ideal is needed in order to guarantee that
$\tau(|x|^p)$ is finite.)
 We refer to \cite{Ne,Ter,Kofa,Yead} for more
on this and the fact that $L_p(N,\tau)$ can be realized as
unbounded operators affiliated to $N$.

\hz

The starting point of Kosaki's \cite{Kos} definition of an
$L_p$-space is a normal faithful state $\phi$ on a von Neumann
algebra $N$. Then $N$ acts on the Hilbert space $L_2(N,\phi)$
obtained by completing $N$ with respect to the norm
 \[ \noo x\rrm_{L_2(N,\phi)} \lel \phi(x^*x)^{\frac12} \pl .\]
The modular operator $\Delta$ is an (unbounded) operator obtained
from the polar decomposition $S=J\Delta^{\frac12}$ of the
antilinear operator $S(x)=x^*$ on $L_2(N,\phi)$ , see
\cite[Section 9.2.]{Kad}. We denote by $\si_t^\phi:N\to N$ the
modular automorphism group defined by $\si_t^\phi(x)\lel
\Delta^{it}x\Delta^{-it}$. Let us recall the standard notation
 \[ x.\phi(y) \lel \phi(xy) \pl .\]
For each $t$ there is a natural map $I_t:N\to N_*$ ($N_*$ the
unique predual of $N$) defined by
 \[ I_t(x)\lel \si_t(x).\phi  \pl . \]
According to \cite[Theorem 2.5.]{Kos}, there is a unique
extension $I_z:N\to N_*$ such that for fixed $x$ the function
$f_x:\{z \p | \p  -1\le Im(z)\le 0\}\to N_*$, $f_x(z) \lel
I_z(x)$ is analytic and satisfies
 \[ f_x(t)(y) \lel \phi(y\si_t^\phi(x)) \quad \mbox{and} \quad
    f_x(-i+t)(y) \lel \phi(\si_t^\phi(x)y) \pl .\]
The density of the algebra of analytic elements shows that for
$0\le \eta \le 1$ the map $I_{-i\eta}$ is injective. By complex
interpolation, the Banach space
 \[ L_p(N,\phi,\eta) \lel [I_{-i\eta}(N),N_*]_{\frac{1}{p}} \]
is defined by specifying $\noo x \rrm_0 \lel \noo
I_{-i\eta}^{-1}(x)\rrm_N$ and $ \noo x\rrm_1 \lel \noo
x\rrm_{N_*}$. We refer to \cite{Ter,Fid} for further information.

\hz

Haagerup's abstract $L_p$ space \cite{Haa,Ter} is defined for
every von Neumann algebra $N$ using the crossed product
$N\rtimes_{\si^w}\rz$ with respect to the modular automorphism
group of a n.s.f. weight $w$. (In our applications, we can assume
$w=\phi$ for a n.f. state.) If $N$ acts faithfully on a Hilbert
space $H$, then the crossed product $N\rtimes_{\si^w}\rz$ is
defined as the von Neumann algebra defined on $L_2(\rz,H)$ and
generated by
 \[ \pi(x)(\xi(t))\lel \si_{-t}^w(\xi(t)) \quad \mbox{and}  \quad
  \la(s)\xi(t)\lel \xi(t-s)  \pl .\]
Then $N\rtimes_{\si^w}\rz$, see \cite{PeTa}, is semifinite and
admits a unique trace $\tau$ such that the dual action
 \[ \theta_s(x) \lel W(s)xW(s)^*\]
satisfies $\tau(\theta_s(x))=e^{-s}\tau(x)$. Here $W(s)$   is
defined by  the phase shift
 \[ W(s)\xi(t) \lel e^{-ist}\xi(t) \pl .\]
The dual action satisfies $\theta_s(\pi(x))=\pi(x)$ and moreover
 \begin{eqnarray}
  \pi(N) &=&   \{ x\in N\rtimes_{\si^w}\rz   \pl |  \pl \theta_s(x)\lel x\pl , \mbox{ for all $s \in \rz$}\} \pl
 \end{eqnarray}
Let us agree to identify $N$ with $\pi(N)$ in the following.
$L_p(N)$ is defined to be the space of unbounded,
$\tau$-measurable  operators
 affiliated to $N\rtimes_{\si^w}\rz$ such that
for all $s\in \rz$
 \[ \theta_s(x) \lel e^{-\frac{s}{p}}x \pl .\]
Note that the intersection $L_p(N)\cap L_q(N)$ is
 $\{0\}$ for  different values  $p\neq q$.
There is a natural isomorphism between $N_*$ and
$L_1(N)$ such that for every normal functional
$\phi\in N_*$ there is a unique $a_\phi \in L_1(N)$
associated satisfying
 \[  \tau(a_{\phi}x)  \lel \phi(\intt_\rz \theta_s(x))
 \pl \]
for all positive $x\in N\rtimes_{\si_t^w}\rz$. The
key point in this construction is the definition of
the trace function $tr:L_1(N)\to \cz$
(corresponding to the integral in the commutative
case) given by
 \[ tr(a_\phi) \lel \phi(1) \pl .\]
Let $\frac1p+\frac{1}{p'} \lel 1$ and $x\in L_{p}(N)$, $y\in
L_{p'}(N)$.  Then we have the trace property
 \[ tr(xy)\lel tr(yx) \pl .\]
The polar decomposition $x=u|x|$ of $x\in L_p(N)$ satisfies $u\in
N$ and
 \[ \noo x\rrm_p \lel tr(|x|^p)^{\frac1p} \pl .\]
$N$ acts as a left and right module on  $L_p(N)$ and more
generally H\"older's inequality
 \begin{eqnarray}
 \noo xy\rrm_r &\le& \noo x\rrm_p \noo y\rrm_q \pl
 \end{eqnarray}
holds whenever $\frac1p +\frac1q=\frac1r$. As for semifinite von
Neumann algebras, there is a positive cone $L_p(N)_+$ in $L_p(N)$
consisting of elements in $L_p(N)$ which are positive as unbounded
operators affiliated to $N\rtimes_{\si^w}\rz$. Following
\cite[Proposition 33, Theorem 32]{Ter}, we deduce for $0\le x\le y
\in L_p(N)$ and $\frac1p+\frac{1}{p'}=1$
 \begin{eqnarray}
  \noo x\rrm_p \lel \sup_{z\in L_{p'}(N)_+, \noo z\rrm_{p'}\le 1} tr(zx)
 &\le&   \sup_{z\in L_{p'}(N)_+, \noo z\rrm_{p'}\le 1} tr(zy)  \lel \noo y\rrm_p \pl .
 \end{eqnarray}
In the sequel, we will often use the following simple
observation.\hhz

\begin{lemma}\label{fact0}  Let $0<p\le \8$ and  $x,y\in L_p(N)$ such that $x^*x\le y^*y$. If
$p$ is the left support projection of $y$, then
$xy^{-1}p$ is a well-defined element in $N$ of norm
less than one.
\end{lemma}\hhz

{\bf Proof:} We note that  $a=xy^{-1}p$ is affiliated with
$N\rtimes_{\si^w}\rz$ and
 \[ p(y^{-1})^*x^*xy^{-1}p \kl p(y^{-1})^*y^*yy^{-1}p \kl p \]
shows that $a$ is a contraction  and in particular
$\tau$-measurable. Moreover, we have
 \[ p \lel \theta_s(p) \lel \theta_s(yy^{-1}p) \lel
 e^{-\frac{s}{p}}y \theta_s(y^{-1}p) \]
and hence
 \begin{eqnarray}
 \theta_s(y^{-1}p) \lel e^{\frac{s}{p}} y^{-1}p \pl .
 \end{eqnarray}
Therefore, the equality
 \for
 \theta_s(a) &=& \theta_s(x)\theta_s(y^{-1}p) \lel
  e^{-\frac{s}{p}} x e^{\frac{s}{p}} y^{-1}p \lel a \pl
  \mel
shows with $(1.1)$ that $a$ is a contraction in $N$.\qed

In the $\si$-finite case, Kosaki's $L_p$-space is isomorphic to
Haagerup's $L_p$ space, see \cite[section 8]{Kos}. Indeed, given
a n. f. state $\phi$ with corresponding density $D$
 in $L_1(N)\cong N_*$, then
 \[ \si_t^\phi(x) \lel D^{it}xD^{-it} \]
and for all $x\in N$
 \for
  \noo I_{-i\eta}(x)\rrm_{[I_{-i\eta}(N),N_*]_{\frac{1}{p}}}
   \lel \noo D^{\frac{\eta}{p}}xD^{\frac{1-\eta}{p}}\rrm_{L_p(N)} \pl .
 \mel
We recall that an element $N$ is analytic, if $t\mapsto
\si_t^{\phi}(x)$ extends to an analytic function with  values in
$N$. The $^*$-closed subalgebra of analytic  elements will be
denoted by $\A$. The following Lemma is probably well-known, see
\cite{Kos}, \cite[Lemma 1.1.]{JX}. We add a short proof for the
convenience of the reader.\hhz

\begin{lemma}\label{Ap0}   Let $0< p<\8$, then
$D^{\frac{1}{2p}}\A_+D^{\frac{1}{2p}}$ is dense in $L_p(N)_+$ and
$ND^{\frac1p}$ is dense in $L_p(N)$ and for $1\le p\le \8$ the map
$J_p:L_p(N)\to L_1(N)$, $J_p(x)=xD^{1-\frac1p}$ is injective.
\end{lemma}\hhz

{\bf Proof:} $D$ is a $\tau$-measurable operator with support
projection $1$. Therefore, $xD^{1-\frac 1p}=0$ implies that
$x=xD^{1-\frac 1p}D^{\frac1p-1}$ is a well-defined
$\tau$-measurable operator and equals $0$. Hence, $J_p$ is
injective for $1\le p\le \8$. Let $\frac1p+\frac{1}{p'}=1$. We
show the density of $D^{\frac{1}{2p}}\A_+D^{\frac{1}{2p}}$ in
$L_p(N)_+$ for $1\le p<\8$. If this is not the case, the
Hahn-Banach theorem implies the existence of $x\in L_p(N)_+$ and
$y\in L_{p'}(N)_{sa}$ such that
 $tr(yD^{\frac{1}{2p}}aD^{\frac{1}{2p}})=tr(D^{\frac{1}{2p}}yD^{\frac{1}{2p}}a)\le
 0$ for all $a\in \A_+$ and $tr(xy)>0$. Using the $\si$-strong
 density of $\A_+$ in $N_+$, see  \cite{PeTa}, we deduce that
$y=D^{-\frac{1}{2p}}D^{\frac{1}{2p}}yD^{\frac{1}{2p}}D^{-\frac{1}{2p}}$
is negative  and hence $tr(xy)\le 0$, a contradiction. Let
$\frac12 \le p\le 1$ and $y\in L_p(N)_+$, then we can approximate
$y^{\frac12}$ by  an element
$x=D^{\frac{1}{4p}}aD^{\frac{1}{4p}}$, $a\in \A_+$ and hence
H\"older's inequality implies that
$x^2-y=x(x-y^{\frac12})+(x-y^{\frac12})y^{\frac12}$ approximates
$y$. Since $a$ is analytic, we observe that $x^2=
 D^{\frac{1}{4p}}aD^{\frac{1}{4p}}D^{\frac{1}{4p}}aD^{\frac{1}{4p}}
 =D^{\frac{1}{2p}}\si_{-\frac{i}{4p}}(a)^*\si_{-\frac{i}{4p}}(a)D^{\frac{1}{2p}}$
is in $D^{\frac{1}{2p}}\A_+D^{\frac{1}{2p}}$. By induction, we
deduce the density of $D^{\frac{1}{2p}}\A_+D^{\frac{1}{2p}}$ in
$L_p(N)_+$ for all $0<p<\8$. Since $\A
D^{\frac1p}=D^{\frac{1}{2p}}\A D^{\frac{1}{2p}}$ and every element
is a linear combination of $4$ positive elements, we obtain  the
assertion.\qed

\section{\bf Hilbert $C^*$-modules and \boldmath$L_p$-spaces\unboldmath}

\newcommand{\Cp}[1]{L_p(#1 ;\ell_2^C)}
\newcommand{\Cc}[2]{L_{#2}(#1 ;\ell_2^C)}
\setcounter{equation}{0} \hhz
 In this section, we will analyze some  Banach
spaces related to the dual version of Doob's inequality and
combine Kasparov's stabilisation theorem for Hilbert
C$^*$-modules with Stein's inequality for non-commutative
martingales  proved by Pisier and Xu \cite{PX}. Let $N\subset M$
be a von Neumann subalgebra and $E:M\to M$ be a normal
conditional expectation onto $N$. $E$ is normal if the dual map
$E^*$ satisfies $E^*(M_*)\subset M_*$ and hence has a predual map
$E_*:M_*\to M_*$ such that $E=(E_*)^*$. In order to simplify the
exposition, we will assume in addition that $\phi:M\to \cz$ is a
normal faithful state satisfying $\phi=\left. \phi\right|_N\circ
E$. Using Kosaki's interpolation spaces \cite[Proposition
4.1]{Kos}, it is very easy to check that $E$ extends to certain
$L_p$-spaces. However, in our context it is more convenient to
work with the Haagerup $L_p$-spaces.\hhz

\begin{lemma}\label{mod1}  Under the previous assumptions let $D$ be the density of $\phi$ in
$L_1(M)$ and $\si_t$ be the modular group, then $E_*(D)\lel D$,
and for all $x\in M$, $y\in M_*$
 \[ E(\si_t(x)) \lel \si_t(E(x)) \quad \mbox{and}  \quad tr(E_*(y)) \lel tr(y)  \pl .\]
\end{lemma}\hhz

{\bf Proof:} The first assertion follows from
 \[  tr(Dx) \lel \phi(x) \lel \phi(E(x))\lel tr(DE(x))\lel tr(E_*(D)x) \]
which is valid for all $x\in M$. For  $E(\si_t(x)) \lel
\si_t(E(x))$  see \cite[Lemme 1.4.3]{C}. In order to prove the
third assertion let $\psi \in M_*$ and $a_\psi$ be the
corresponding density then
 \for
   tr(E_*(a_\psi)) &=&  \psi\circ E(1) \lel \psi(E(1)) \lel \psi(1) \lel tr(a_\psi)
  \pl .\\[-1.65cm]
   \mel \qed

We will  need several approximation results.\hhz

\begin{lemma}\label{Ap1}  Let $\A$ be a $^*$-closed (not necessarily norm closed) but
$\si$-strongly dense subalgebra of $M$. For $0< p<\8$ the space
$\A D^{\frac{1}{p}}$ is norm dense in $L_p(N)$.
\end{lemma}\hhz

{\bf Proof:} According to Lemma \ref{Ap0}, $M D^{\frac{1}{p}}$ is
norm dense in $L_p(M)$. According to Kaplansky's density theorem
\cite[Theorem II.4.8]{TAK}, every element in the unit ball of $M$
is the strong$^*$-limit of elements in the unit ball of $\A$. An
application of the following Lemma \ref{Ap2} yields the
assertion.\qed

\begin{lemma}\label{Ap2}  Let $0<p<\8$, $M$ be a von Neumann algebra,
$a\in L_1(M)_+$ ,$(x_\al)_\al\subset N$ be a bounded net and $x\in
N$. If $x_\al a^{\frac1p}$ converges to $xa^{\frac1p}$ in norm for
some $1\le p<\8$, then this is true for all $0<p<\8$. In
particular, if $(x_\al)_\al\subset N$ converges to $x$ strongly,
then for all $b\in L_p(N)$, $x_\al b$ converges to $xb$ in norm.
\end{lemma}\hhz

{\bf Proof:} We  consider the set $I \lel  \left \{ 0<p<\8 \mitt
\lim_\al x_\al a^{\frac 1p} \lel xa^{\frac 1p}  \right \}$. Here,
we refer to convergence in norm. Let us first observe that
$0<q<p\in I$ implies $q\in I$. Indeed,  according to H\"{o}lder's
inequality $(1.2)$,   the linear map $M_r: L_p(N)\to L_q(N)$
defined by $M_l(x) \lel xa^{\frac 1q -\frac 1p}$ is continuous
and therefore, the convergence  of $x_\al a^{\frac 1p}$ implies
the convergence of $x_\al a^{\frac1q}=M_r(x_\al a^{\frac 1p})$.
Let us now show that $1\le p\in I$ implies $2p\in I$, i.e.
 \[ \noo (x_\al -x)a^{\frac{1}{2p}}\rrm_{2p} \lel 0  .\]
We note  that
\[ \noo (x_\al -x)a^{\frac{1}{2p}}\rrm_{2p}^2 \lel
 \noo a^{\frac{1}{2p}}(x_\al^* -x^*)(x_\al-x) a^{\frac{1}{2p}}
 \rrm_{p} \]
and define the analytic function $f:\{0\le Im(z)\le 1\}\to L_p(N)$
given by
 \[ f_\al (z) \lel \exp((z-\frac 12)^2) \pl   a^{\frac{z}{p}}
 (x_\al^*
 -x^*)(x_\al-x) a^{\frac{1-z}{p}} \pl .\]
Since $a^{it}xa^{-it}\in L_p(N)$ for every  $x\in L_p(N)$ this is 
well defined and  analytic in the interior. For $z=1+it$, we get
 \for
 \noo f_\al(1+it)\rrm_p &\le& \exp( \frac 14) \pl
  \noo a^{\frac{1+it}{p}} (x_\al^* -x^*)(x_\al-x)
  a^{\frac{-it}{p}}\rrm_p \\
  &\le& \exp( \frac 14) \pl
  \noo  a^{\frac{1}{p}}(x_\al^* -x^*)\rrm_p \pl \noo
  x_\al-x\rrm_\8
  \\
  &\kl &  \exp( \frac 14) \pl \noo (x_\al -x)a^{\frac{1}{p}} \rrm_p \pl \sup_\al \noo
  x_\al-x\rrm_\8
 \pl .
  \mel
Similarly,
  \for
 \noo f_\al(it)\rrm_p &\le& \exp( \frac 14) \pl
  \noo a^{\frac{it}{p}} (x_\al^* -x^*)(x_\al-x)
  a^{\frac{1-it}{p}}\rrm_p \\
  &\le& \exp( \frac 14) \pl \sup_\al \noo x_\al^* -x^*\rrm_\8
  \noo (x_\al-x)  a^{\frac{1}{p}}\rrm_p  \pl .
  \mel
Therefore, $p\in I$ implies with  the three line Lemma
 \for
 \lefteqn{  \lim_\al \noo   a^{\frac{1}{2p}}(x_\al^*-x^*)(x_\al-x)
 a^{\frac{1}{2p}}
 \rrm_{p} \lel  \lim_\al \noo f_\al(\frac 12) \rrm_p  }\\
  & & \kl \lim_\al \sup_t \max\{ \noo f_\al(it)  \rrm_p, \noo
  f_\al(1+it)\rrm_p\} \\
  & & \kl  \exp(\frac 14)\pl  \sup_\al \noo x_\al-x\rrm_\8 \pl
  \lim_\al \noo
 (x_\al-x) a^{\frac{1}{p}} \rrm_p   \lel 0 \pl .
  \mel
Hence, if $I\cap [1,\8)$ is non empty, then $I=(0,\8)$ and the
first assertion is proved.  Let $b\in L_p(N)$.  We use the polar
decomposition of $b^*$ to find a partial isometry $u$ such that
$b=|b^*|u$. Then $a=|b^*|^p$ is in $L_1(N)_+$.  If $x_\al$
converges to $x$ strongly then  $x_\al a^{\frac12}$ converges in
norm to $xa^{\frac12}$. This means $2\in I$ and hence $I=(0,\8)$.
In particular, $p\in I$ and hence $x_\al|b^*|=x_\al a^{\frac1p}$
converges to $x|b^*|$. This immediately implies that $x_\al b=
x_\al|b^*|u$ converges to $xb=x|b^*|u$.\qed

In this paper, we will frequently use the space $\Cp M \subset
L_p(B(\ell_2)\bar{\ten} M)$ of column matrices with values in
$L_p(M)$.\hhz

\begin{cor} \label{Ap3} Let $0<p<\8$ and $\A$ $\si$-strongly dense $^*$-closed subalgebra of $M$, then
the space of finite sequences $(a_nD^{\frac1p})_n$ such that
$a_n\in \A$ is dense in $\Cp M $.
\end{cor}\hhz

{\bf Proof:} Let $p_n$ be the orthorgonal projection onto the
first $n$ unit vectors in $\ell_2$, then $p_n\ten 1$ converges to
$1\ten 1$ in the strong$^*$ topology. Hence, the space of finite
sequences is dense in $\Cp M $. We can apply Lemma \ref{Ap1} for
each coordinate to obtain the assertion. \qed

In the next step we refer to \cite[Proposition 2.3]{JX} to
understand how a state preserving conditional expectation extends
to a contraction $E_p$ on $L_p(M)$ and in fact $E_*=E_1$. We will
later often use these facts without further reference and drop
the index $p$ or $_*$  because it is easily determined by the
context. \hhz

\begin{prop}\label{baby}  Let $1\le p \le \8$, then there is
a contraction $E_p:L_p(M)\to L_p(M)$ such that for all $x\in M$
and $0\le \theta \le 1$
 \[ E_p(D^{\frac{1-\theta}{p}}xD^{\frac{\theta}{p}})\lel
  D^{\frac{1-\theta}{p}}E(x)D^{\frac{\theta}{p}} \]
$E_p$ satisfies $E_p(x^*)=E_p(x)$ and $E_p(x)$ is positive for
all positive  $x$.  Moreover, if $\frac 1s=\frac 1p+\frac
1q+\frac 1r\le 1$, $a\in L_p(N)$, $b\in L_q(N)$ and $x\in
L_r(M)$, then
 \[  E_s(axb) \lel a E_r(x) b \pl .\]
\end{prop}\hhz

Conditional expectations (or, more generally, positive
operator-valued weights) are closely connected to Hilbert
C$^*$-modules. An excellent reference for the few facts and the
notation we need  in this paper is Lance's book \cite{LA}. We
recall that a Hilbert C$^*$-module $\M$ over a $C^*$-algebra $\A$
is a right $\A$-module $\M$ together with an $\A$-valued
sesquilinear form $\langle \cdot,\cdot\rangle:\M\times \M\to \A$
such that for all $a\in \A$ and $x,y\in \M$
 \for
 \langle x,ya \rangle &=&  \langle x,y \rangle a \pl ,\\
 \langle x,y \rangle^* &=& \langle y,x \rangle \pl .
 \mel
The norm in $\M$ is given by $\noo x\rrm_\M=\noo \langle
x,x\rangle\rrm_\A^{\frac12}$.  We should note that the
sesquilinear form is assumed to be linear in its second and
antilinear in its first component. The standard example of a
Hilbert $C^*$-module over a $C^*$-algebra $\A$ is the space
$H_\A$, the closure of finite sequences $(x_n)\subset \A$ with
respect to the norm
 \[ \noo (x_n)\rrm \lel \noo \summ_n x_n^*x_n\rrm_\A^{\frac12} \pl
 .\]
We can identify $H_\A$ with the space of column matrices $H_\A
\subset \K\ten \A$. In the context of von Neumann algebras $\A=N$,
we derive a similar example using the space of column matrices
$C(N) \subset B(\ell_2) \bar{\ten}N$ which consists of sequences
$(x_n)$ such that $\sum_n x_n^*x_n$ converges in the $\si$-weak
operator topology. Then the sesquilinear form
 \[ \langle (x_n), (y_n)\rangle \lel (x_1^*,x_2^*,\cdots ) (y_1,y_n,\cdots)^{\perp} \lel \summ_n x_n^*y_n \]
converges in the $\si$-weak operator topology and satisfies the
axioms of a Hilbert $C^*$-module. Let us indicate how this concept
applies to non-commutative $L_p$-spaces. Indeed, $Ses :\Cp N\times
\Cc N q \to L_r(N)$ defined by
 \[ Ses ((x_n),(y_n)) \lel \summ_n x_n^*y_n \quad, \quad
 \mbox{where}\quad
 \frac1r\lel \frac1p+\frac1q \]
is still a sesquilinear form. We will omit $Ses$ and simply write
$(x_n)^*(y_n)$ to remind the natural matrix multiplication.
H\"older's inequality immediately implies
 \begin{eqnarray}
  \noo \summ_n x_n^*y_n\rrm_r &\le&  \noo (\summ_n x_n^*x_n)^{\frac12}
 \rrm_p \noo (\summ_n y_n^*y_n)^{\frac12}
 \rrm_q \pl
 \end{eqnarray}
whenever $\frac1r=\frac1p+\frac 1q$. Our aim is to extend this
type of inequalities to the sesquilinear form
   \[ \langle x,y\rangle \lel E(x^*y) \]
given by a conditional expectation $E:M\to N$. We denote by
$L_\8(M,E)$ the completion of $M$ with  respect to the norm
 \[ \noo x\rrm_{L_\8(M,E)}  \lel \noo E(x^*x)\rrm_\8^{\frac12} \pl .\]
In analogy with the case $p=\8$, we introduce the following
notation:\hhz

\begin{defi} Let $0< p \le \8$ and $x=aD^{\frac1p} \in MD^{\frac1p}$, then
 \[ \noo x\rrm_{L_p(M,E)} \lel \noo D^{\frac1p}E(a^*a)D^{\frac1p}\rrm_{\frac{p}{2}}^{\frac 12} \pl .\]
The completion with respect to this 'norm' is denoted by
$L_p(M,E)$.
\end{defi}\hhz

\begin{rem} \label{prec} {\rm Let us note that for $p \ge 2$, the precaution $x\in
MD^{\frac1p}$ is unnecessary because for all $x\in L_p(N)$ the
conditional expectation  $E(x^*x)$ is well-defined and
 \[ \noo x\rrm_{L_p(N,E)}\lel \noo E(x^*x)\rrm_{\frac p2}^{\frac12} \kl \noo x^*x\rrm_{\frac p2}^{\frac12}\lel \noo x\rrm_p \pl .\]
By norm density of $MD^{\frac1p}$ in $L_p(M)$, we obtain the same
closure. However, for $p<2$ we no longer dispose of the continuity
of $E$ on $L_{\frac{p}{2}}(N)$ (even in the commutative case).
This justifies our slightly artificial definition.}
\end{rem}\hhz

Our next step will be to prove the triangle inequality for $p\ge
1$  using Kasparov's stabilisation result.\hhz

\begin{prop} \label{CS1} Let $N\subset M$, $E$ and $\phi$ as
above and assume in  addition that $M_*$ is
separable. For all $0<p\le \8$, there exists an
isometric right $N$-module map $u_p:L_p(M,E)\to \Cp
N$.  Moreover,
 \begin{enumerate}
 \item[i)] For all $x\in L_p(M,E)$ and $y\in L_q(M,E)$
  \[ u_p(x)^*u_q(y)\lel E(x^*y) \pl .\]
 \item[ii)] Let $\frac1r=\frac1q+\frac1p$. For all $x\in L_p(M,E)$ and $y\in L_q(M,E)$
  \[ \noo E(x^*y)\rrm_r \kl \noo x\rrm_{L_p(M,E)} \noo
  y\rrm_{L_q(M,E)}\pl .\]
 \item[iii)] If $0<p<\8$ there exists a contractive projection $Q_p$ onto the image of $u_p$ such that for all
 $z\in \Cp N$
 \[   Q_p(z)^*Q_p(z)\kl z^*z  \pl .\]
 \item[iv)] If $1<p\,, \, p'<\8$ and $\frac 1p+\frac{1}{p'}=1$, then
  \[ Q_p^*\lel Q_{p'} \pl .\]
 \end{enumerate}
In particular, $L_p(M,E)$ is a normed space for $1\le p<\8$ and
$p$-normed for $0<p\le 1$.
\end{prop}\hhz

{\bf Proof:} Let $X\subset M$ be a $\si$-strongly dense separable
$C^*$-algebra and let $F$ be the Hilbert $C^*$-module generated
by $N$ and $X$. Then $F$ is a countably generated Hilbert
$N$-module and according to \cite[Theorem 6.2.]{LA} there exists
a unitary $w: F\oplus H_N\to H_N$ such that
 \[ w(f,x)^*w(g,y) \lel \langle f,g\rangle + x^*y \pl .\]
Let $\tilde{Q}:F\oplus H_N\to F$ be the projection onto the first
coordinate. This projection carries over to
$Q=w\tilde{Q}w^{-1}:H_N\to w(F)$. We  define $u:F\to H_N$ to be
the restriction of $w$ to the first component. Clearly, $u$
preserves the $N$-valued sesquilinear form. For $0<p \le \8$,  we
define the map $u_p:FD^{\frac1p} \to \Cp N$ by
 \[ u_p(aD^{\frac1p} ) \lel u(a)D^{\frac1p} \pl .\]
Then for all $a,b\in F$
 \[  D^{\frac1p}E(a^*b)
 D^{\frac1q} \lel D^{\frac1p}u(a)^*u(b)D^{\frac1q}  \lel
 u_p(aD^{\frac1p})^*u_q(bD^{\frac1q}) \pl . \]
This justifies $i)$ for the subsets $FD^{\frac 1p}$,
$FD^{\frac1q}$, respectively. In particular, we obtain
  \[ \noo a D^{\frac 1p}\rrm_{L_p(M,E)}^2
   \lel \noo D^{\frac 1p}E(a^*a)D^{\frac 1p}\rrm_{\frac p2}
   \lel \noo u_p(aD^{\frac 1p})^*u_p(aD^{\frac 1p})\rrm_{\frac p2}
    \lel \noo u_p(xD^{\frac 1p})\rrm_p^2 \pl \]
which implies that $u_p$ is isometric when restricted to
$FD^{\frac1p}$. According to the following Lemma \ref{Ap4} and
Proposition \ref{p8}, $FD^{\frac1p}$ is dense and hence $i)$
follows. Assertion $ii)$ is an immediate consequence of $(2.1)$.
Indeed,
 \for
 \noo D^{\frac1p}E(a^*b) D^{\frac1q}\rrm_r &=&
 \noo u_p(aD^{\frac1p})^*u_q(bD^{\frac1q}) \rrm_r \kl
 \noo u_p(aD^{\frac1p})\rrm_p  \noo u_q(bD^{\frac1q}) \rrm_q \\
 &=&
  \noo aD^{\frac1p}\rrm_{L_p(M,E)} \pl  \noo
  bD^{\frac1p}\rrm_{L_q(M,E)}\pl .
 \mel
Since $\Cp N$ is a Banach space, respectively a complete
$p$-normed space, and $u_p$ is isometric, we deduce the last
assertion. The projection $Q_p :\Cp N \to \Cp N$ is densely
defined for a sequence $(z_nD^{\frac 1p})$ such that $z_n=0$ for
$n\ge n_0$ by
 \[  Q_p[(z_nD^{\frac 1p})] \lel Q[(z_n)] D^{\frac 1p} \pl .\]
(Here and in the following we will use the brackets $[ \pl ]$ to
indicate that $Q$ is applied to the sequence $(z_n)$.) Since  $w$
is a unitary in the sense of Hilbert $N$-modules, we deduce for
$w^{-1}[(z_n)]=f + (a_n)$
 \for
  Q[(z_n)]^*Q[(z_n)] &=& \langle \tilde{Q}(f+(a_n)), \tilde{Q}(f+(a_n)) \rangle
   \lel \langle f,f\rangle
   \kl \langle f,f\rangle + \summ_n a_n^*a_n \\
   &=& \langle f + (a_n), f + (a_n)\rangle \lel
      \langle w(f + (a_n)), w(f + (a_n)) \rangle  \lel  \summ_n z_n^*z_n
      \pl.
  \mel
In particular,
 \for
  \noo Q_p[(z_nD^{\frac 1p})]\rrm_p^2  &=&  \noo D^{\frac 1p} Q[(z_n)]^*Q[(z_n)]D^{\frac 1p}
 \rrm_{\frac{p}{2}}  \kl
  \noo D^{\frac 1p} \summ_n z_n^*z_n D^{\frac 1p}
 \rrm_{\frac{p}{2}}\\
 &=&
 \noo (z_nD^{\frac 1p})\rrm_{\Cp N}^2 \pl .
 \mel
Hence, $Q_p$ extends by density, see Corollary \ref{Ap3}, to a
contraction on $\Cp N$. This proves $iii)$. In order to prove
$iv)$ it suffices again by density to consider finite sequences
$z=(z_nD^{\frac1p})$ and
$\tilde{z}=(\tilde{z}_nD^{\frac{1}{p'}})$ with $z_n,
\tilde{z_n}\in F$. Using again $w^{-1}[(z_n)]=f + (a_n)$,
$w^{-1}[(\tilde{z}_n)]=\tilde{f} + (\tilde{a}_n)$, we deduce from
the fact that $w$ is a unitary module map
  \for
 tr(Q_p(z)^*\tilde{z}) &=& tr(D^{\frac1p}Q[(z_n)]^*(\tilde{z}_n)D^{\frac{1}{p'}})  \\
  &=&   \phi( \langle Q[(z_n)], (\tilde{z}_n)\rangle) \\
  &=& \phi( \langle \tilde{Q}(f + (y_n)), \tilde{f}+
  (\tilde{a}_n)\rangle) \\
  &=&
 \phi( E(f^*\tilde{f})) \\
  &=&  \phi( \langle f + (a_n), \tilde{Q}(\tilde{f}+
  (\tilde{a}_n)) \pl \rangle)  \\
  &=&
  tr(z^*Q_{p'}(\tilde{z})) \pl .\\[-1.675cm]
  \mel\qed

Let us briefly clarify how the unitary module map  $w: F\oplus
H_N\to H_N$ extends to an isometry onto $C(N)$. We recall that
$M$ acts naturally on $L_2(M)$ and that $p_E=E_2:L_2(M)\to
L_2(M)$ is the orthogonal projection onto $L_2(N)$. Hence, we
obtain an isometric isomorphism
  \[ L_\8(M,E) \lel Mp_E  \subset B(L_2(N),L_2(M)) \subset B(L_2(M))  \pl .\]
Note that in general $p_E$ is not in $M$. Via this inclusion, we
dispose of all the relevant locally convex operator topologies on
$L_\8(M,E)$ and on
 \[ L_\8(M,E) \oplus C(N) \subset B(L_2(N), L_2(M)\oplus
 \ell_2(L_2(N))) \pl .\]

We denote by $L_\8^{st}(M,E)$ the closure with respect to the
$\si$-strong topology. The following example due to O. Ramcke
shows that this might be different from $L_\8(M,E)$. Let
$M=L_\8([0,1]^2)$ and $E(f)(t,s)=\int_0^1 f(t,r) dr$. Then, the
norm on $L_\8(M,E)$ is the norm in $L_\8(L_2)$. Given disjoint
sets $I_k$ of measure $2^{-k}$ the function
 \[ f \lel \summ_k 1_{I_k} \ten 2^{\frac{k}{2}} 1_{I_k} \]
yields an element in $\si$-strong closure not belonging to
$L_\8(M,E)$. The next proposition shows that $L_\8^{st}(M,E)$ is
isomorphic to a complemented module in $C(N)$. In particular, we
will find a sequence $(x_j)\subset L_\8^{st}(M,E)$ and completely
contractive module maps $u_j:M\to N$ such that for all $x\in M$
 \[ x\lel \summ_j x_ju_j(x) \]
For inclusions with finite index, a finite sum of this form
suffices and more can be said about the coefficients
$E(x_i^*x_i)$, see \cite{Bi}. In any case
 \[ 1 \lel E(1) \lel \summ_j u_j(1)^*E(x_j^*x_j)u_j(1) \]
shows that we are dealing with a 'partition of unity'.\hhz

\begin{prop} \label{p8} Under the above assumptions,  $w$ extends to an isometric isomorphism between
the space $L_\8^{st}(M,E) \oplus C(N)$ and $C(N)$. Moreover, for
every element $x$ in the unit ball of $L_\8(M,E)$ there exists a
net $x_\al$ in the unit ball of $L_\8(M,E) \cap F$ converging to
$x$ in the strong$^*$ topology and $Q$ extends to a projection
onto the image of $L_\8^{st}(M,E)$.
\end{prop}\hhz

{\bf Proof:} We note that for $h\in L_2(N)$ and $x\in F\oplus
C(N)$
  \[ \noo w(x)(h)\rrm^2 \lel (w(x)^*w(x)(h),h) \lel (x^*x(h),h)
  \lel \noo x(h)\rrm_2^2 \pl .\]
Hence, $w$ preserves the strong and $\si$-strong
topology. Therefore, we obtain a natural extension,
also denoted by $w$, to the $\si$-strong closure
$L_\8^{st}(M,E)\oplus C(N)$  of $F\oplus H_N$ with
values in $C(N)$. To show that $w$, $w^{-1}$  remain
contractions, it suffices to show that every
element in the unit ball of $L_\8^{st}(M,E)\oplus
C(N)$,  $C(N)$,  can be approximated with respect
to the strong topology by elements in the unit ball
of $F\oplus H_N$, $H_N$,  respectively. This follows
immediately from Kaplansky's density theorem
\cite[Theorem II.4.8]{TAK} for $C(N)$. To obtain
the assertion for the $\si$-strong closure of
$Mp_E$, we note that $F$ is $\si$-strong$^*$ dense
in $M$, hence $Fp_E$ is $\si$-strong$^*$ dense in
$L_\8^{st}(M,E)$. Then, we follow the proof of
Kaplansky's density theorem,  see  \cite[Theorem
II.4.8]{TAK}, and note that the  function $f:F\to
B(L_2(N))$ defined there by
 \[ f(xp_E)\lel 2xp_E(1+(p_Ex^*xp_E))^{-1} \lel
   2xp_E(1+(E(x^*x))^{-1} \in FN\subset F  \pl \]
satisfies $f(F)\subset F$. Hence, the  proof of
Kaplansky's density theorem applies  and   therefore
every element in the unit ball of
 \[  L_\8^{st}(M,E) \lel \overline{Mp_E}^{\si-\mbox{\scriptsize strong}}
  \lel \overline{Mp_E}^{\si-\mbox{\scriptsize strong}^*}  \]
(\cite[Theorem 2.6]{TAK}) can be approximated by a net $f(x_\al)$
in the unit ball of $F$. Since $Q=w^{-1}\tilde{Q}w$ corresponds
to the projection onto the first component in $L_\8^{st}(M,E)
\oplus C(N)$, the last assertion follows easily. \qed

\begin{lemma}\label{Ap4} Under the previous assumptions $FD^{\frac1p}$ is dense in $L_p(M,E)$.
The Cauchy-Schwarz inequality \ref{CS1} $ii)$ and $i)$ also holds
for $p=\8$ or $q=\8$.
\end{lemma}\hhz

{\bf Proof:} Let us assume $2\le p <\8$ first. Using Kaplansky's
density theorem \cite[Theorem II.4.8]{TAK}, the unit ball of $F$
is strong$^*$-dense in the unit ball of $M$. Using Lemma
\ref{Ap2}, we deduce that $FD^{\frac1p}$ is dense in $L_p(M)$. As
observed in Remark \ref{prec}, the inclusion $L_p(M)\subset
L_p(M,E)$ is contractive and dense, hence  the assertion follows.
For $p<2$, we observe that the inclusion $i:L_2(M,E)\subset
L_p(M,E)$, $i(aD^{\frac12})= aD^{\frac1p}$ is contractive using
H\"older's inequality
 \[ \noo aD^{\frac 1p}\rrm_p \lel
 \noo D^{\frac 1p}E(a^*a)D^{\frac 1p}\rrm_{\frac{p}{2}}^{\frac12}
 \kl  \noo D^{\frac 12}E(a^*a)D^{\frac 12}\rrm_1^{\frac12}
 \lel \noo aD^{\frac12}\rrm_{L_2(M,E)} \pl . \]
By definition of $L_p(M,E)$, the image of $i$ is dense. Hence,
$FD^{\frac1p}=i(FD^{\frac12})$ is dense in $L_p(M,E)$. To prove
the Cauchy-Schwarz inequality if $p=r<\8$ and $q=\8$, we fix
$x=bD^{\frac1p}$ and $y\in L_\8^{st}(M,E)$ of norm less than one.
Then there exists a net $a_\al p_E$ in the unit ball of $F$ such
that $a_\al$ converges to $y$ with respect to the strong$^*$
topology. Then the strong$^*$ convergence of $b^*a_\al p_E$
implies the strong$^*$ convergence of $E(b^*a_\al)=p_Eb^*a_\al
p_E$. Therefore the norm convergence of $D^{\frac1p}E(b^*a_\al)$
to $D^{\frac1p}E(b^*y)$ follows from Lemma \ref{Ap2}. In
particular according to  Proposition \ref{CS} $ii)$,
 \for
 \noo D^{\frac1p}E(b^*y)\rrm_{p} &=&
   \lim_\al \noo D^{\frac1p}E(b^*a_\al)\rrm_{p}
   \kl  \lim_\al \noo D^{\frac1p}E(b^*b)D^{\frac1p}\rrm_{\frac p2}^{\frac12}
    \noo E(a_\al^*a_\al)\rrm_\8^{\frac12} \\
  &\le&  \noo bD^{\frac1p}\rrm_{L_p(M,E)}    \pl .
  \mel
Thus, by density, every norm one element in $L_\8^{st}(M,E)$
induces a contractive 'multiplier' on $L_p(M,E)$. The case $p=\8$
and $q<\8$ is similar. The case $p=\8$ and $q=\8$ is classical and
follows from the fact that $E:M\to B(L_2(M))$ is completely
positive and therefore admits a dilation $E(x)=v\pi(x)v^*$ for a
contraction $v$ and a $^*$-representation $\pi$. Hence, we obtain
\begin{samepage}  \for
  \noo E(x^*y)\rrm_{B(H)} &=&  \noo v\pi(x^*y)v^*\rrm \kl
  \noo v\pi(x^*)\pi(x)v^*\rrm^{\frac12}  \noo
  v\pi(y^*)\pi(y)v^*\rrm^{\frac12}\\
  &=& \noo E(x^*x)\rrm^{\frac12} \pl \noo E(y^*y)\rrm^{\frac12}
  \pl .\\[-1.7cm]
  \mel \qed \end{samepage}

Now we turn our attention to the duality between $L_p(N,E)$ and
$L_{p'}(N,E)$. Let us point out that we use the antilinear duality
bracket
 \[ (x,y) \lel tr(x^*y) \]
between $L_p(N;\ell_2^C)$ and $L_{p'}(N;\ell_2^C)$. The following
Lemma is a standard application of the Hahn-Banach theorem.\hhz

\begin{lemma} \label{dual0} Let $1<p<\8$ and $\frac 1p+\frac{
1}{p'}$, $X_p\subset \Cp M$ and $X_{p'}\subset \Cc M {p'}$ be
subspaces and $Q_p: \Cp M\to \Cp M $ be a projection onto $X_p$
such that $Q_p^*$ is a projection onto $X_{p'}$, then
 \[ X_p^* \lel X_{p'} \pl .\]
and for every dense subset $X\subset X_{p'}$ and $x\in X_p$
 \[ \noo x\rrm \kl \noo Q\rrm \pl \sup \left\{ tr(y^*x) \mitt y\in X\p
 ,\p \noo y\rrm_{X_{p'}} <1  \right\} \pl .\]
\end{lemma}\hhz

{\bf Proof:} Since $(\Cp M)^*=\Cc M {p'}$, the antilinear map
$\iota:X_{p'}\to  X_p^*$
 \[ \iota(y)(x) \lel tr(y^*x) \]
is obviously  contractive.  Let $f:X\to \cz$ be a norm one
functional, then we can apply the Hahn-Banach extension theorem
and $(\Cp M)^*=\Cc M {p'}$ to obtain $z\in \Cc M {p'}$ such that
   \[ f(x) \lel tr(z^*x) \lel ( z, x)  \pl  .\]
Clearly,
 \for
 (z, x)   &=& ( z, Q_p(x)) \lel (Q_p^*(z), x)
 \lel  tr( Q_p^*(z)^*x) \pl.
 \mel
(Defined  in this way $Q_p^*$ is linear.) Since, $Q_p^*(z)\in
X_{p'}$ we deduce that $\iota$ is surjective. The last formula
follows from the Hahn-Banach theorem, i.e.
 \for
  \noo x\rrm &=& \sup_{\noo f\rrm\le 1} |f(x)| \lel \sup_{\noo
  z\rrm \le 1} |tr(z^*x)| \lel
  \sup_{\noo  z\rrm \le 1} |tr(Q_p^*(z)^*x)| \\
  &\le& \noo Q_p\rrm   \sup_{\noo  y\rrm \le 1, y\in X_{p'} } |tr(y^*x)|
 \mel
The supremum is unchanged if  restricted to a dense subset. \qed

\begin{cor}
Let $1<p, p'<\8$, $\frac1p+\frac{1}{p'}=1$  and
$M_*$ be separable,  then $L_p(M,E)^*=L_{p'}(M,E)$
holds isometrically.
\end{cor}\hhz

In the last part of this section, we investigate the space
$L_p(M,(E_n);\ell_2^C)$, a generalization of $L_p(M,E)$,  which
is  important for the dual version of Doob's inequality. We feel
that the space $L_p(M,E)$ is easier to understand and more
directly    connected  to Hilbert $C^*$-modules.  Let $(E_n)$ be
a sequence of conditional expectations $E_n:M\to M$ onto von
Neumann subalgebras $N_n$ such that $\phi \circ E_n =\phi$ for
all $n\in \nz$. Pisier, Xu proved a non-commutative version of
Stein's inequality. To formulate the version we need here, we
consider the subspace $L_p^{cond}(M;\ell_2^C) \subset
L_p(M;\ell_2^C(\nz^2))$ of double indexed sequences $(x_{nk})$
such that $x_{nk} \in L_p(N_n)$ for all $k\in \nz$. We refer to
\cite{PX,JX} for the proof of the following theorem.\hhz

\begin{theorem}[Stein's inequality]\label{Stein}  Let
 $1< p\,,\, p'<\8$ such that  $\frac1p+\frac{1}{p'}=1$ and $(E_n)$ be a
sequence of conditional expectations such that
$E_nE_m=E_{\min\{n,m\}}$. Then the linear map
$ST_p:L_p(M;\ell_2^C(\nz^2))\to L_p^{cond}(M;\ell_2^C)$ defined
by $ST_p[(x_{nk})]=(E_n(x_{nk}))$ is a bounded projection onto
$L_p^{cond}(M;\ell_2^C)$ and satisfies $ST_p^*=ST_{p'}$. In
particular,
 \[ L_p^{cond}(M;\ell_2^C)^*  \lel  L_{p'}^{cond}(M;\ell_2^C) \]
with equivalent norms  depending only on $p$.
\end{theorem}\hhz

\begin{defi} Let $0<p\le \8$, the  space $L_p(M,(E_n);\ell_2^C)$ is the completion
of the space of sequences $(a_nD^{\frac1p})$, $a_n\in M$  with
respect to the norm
 \[ \noo (a_nD^{\frac1p})\rrm_{L_p(M,(E_n);\ell_2^C)}
  \lel \noo D^{\frac1p}\summ_n E_n(a_n^*a_n)D^{\frac1p}\rrm_{\frac p 2}^{\frac12}\pl.\]
\end{defi}\hhz

\begin{prop}\label{CS}
Let $0< p \le\8$ and in addition $M_*$ separable,  then there is
an isometric embedding $u_p:L_p(M,(E_n);\ell_2^C) \to
L_p^{cond}(M;\ell_2^C)$ and for $p<\8$ a norm one projection $R_p$
onto the image of $u_p$. Moreover,
 \begin{enumerate}
  \item[i)] if  $\frac1p+\frac1q=\frac1r$, then for all $(x_n)\in L_p(M,(E_n);\ell_2^C)$
  and $(y_n)\in L_q(M,(E_n);\ell_2^C)$
   \[ u_p(x_n)^*u_q(y_n) \lel \summ_n E_n(x_n^*y_n) \]
 and
 \[ \noo \summ_n E_n(y_n^*x_n)\rrm_r \kl \noo
 (y_n)\rrm_{L_q(M,(E_n);\ell_2^C)} \pl
 \noo (x_n)\rrm_{L_p(M,(E_n);\ell_2^C)} \pl .\]
 \item[ii)] if $1<p,p'<\8$ such that  $\frac1p+\frac{1}{p'}=1$  and
 $ST_p$, $ST_{p'}$ denote the projection onto $L_p^{cond}(M;\ell_2^C)$,
 $L_{p'}^{cond}(M;\ell_2^C)$, respectively, then
  \[ (R_pST_p)^* \lel R_{p'}ST_{p'} \pl .\]
 \end{enumerate}
\end{prop}\hhz

{\bf Proof:} The proof of the Cauchy-Schwarz inequality for
$p=q=\8$ and the isometric embedding  is similar to the last part
of the proof in Lemma \ref{Ap4}. Since we will not need it, we
omit the details. For each $\nen$, we fix the  isometric
isomorphism $u_p^n:L_p(M,E_n) \to L_p(N_n;\ell_2^C)$  and define
$u_p: L_p(M,(E_n);\ell_2^C) \to L_p^{cond}(M;\ell_2^C)$ by
 \[ u_p[(x_n)] \lel ( u_p^n(x_n)_{kn} )\pl ,  \]
i.e. we apply $u_p^n$ to $x_n$ and obtain a double indexed
sequence. For finite sequences, we apply Proposition \ref{CS1}
$i)$ and obtain
 \begin{eqnarray}
  u_p[(x_n)]^*u_q[(y_n)] &=&  \summ_n u_p^n(x_n)^*u_q^n(x_n) \lel  \summ_n E_n(x_n^*y_n) \pl .
 \end{eqnarray}
In particular, $u_p$ is isometric when restricted to finite
sequences. For $p<\8$ an easy Cauchy sequence argument implies
that $u_p$ extends isometrically to $L_p(M,(E_n);\ell_2^C)$.
Similarly as in the proof of Proposition \ref{CS1}, we deduce from
$(2.2)$ the Cauchy-Schwarz inequality. If $Q_p^n:
L_p(N_n;\ell_2^C)\to L_p(N_n;\ell_2^C)$ denotes the projection
onto the image of $u_p^n$, then
 \[ R_p[(x_{nk})] \lel  (Q_p^n[ (x_{nk})_{k}]_{kn}) \]
is certainly well-defined for finite sequences. However, we deduce
from Proposition \ref{CS1} $iii)$
  \[ R_p[(x_{nk})]^*R_p[(x_{nk})]\lel \summ_{n} \summ_k
  Q_p^n[(x_{nk})_{k}]_k^*Q_p^n[(x_{nk})_{k}]_k \kl
  \summ_n \summ_k x_{nk}^*x_{nk} \pl  \]
and therefore $R_p$ extends to a contraction for all $p<\8$.
Finally, we observe that Proposition \ref{CS1} $iii)$ implies
assertion $ii)$:
 \for
  tr( Q_p[(E_n(x_{nk})) ]^*  (y_{nk}) )
  &=& \summ_{n} \summ_k  tr(Q_p^n[(E_n(x_{nk}))_k  ]_k^* y_{nk}) \\
  &=& \summ_{n} \summ_k  tr(Q_p^n[(E_n(x_{nk}))_{k}]_k^* E_n(y_{nk})) \\
  &=& \summ_{n} \summ_k  tr(E_n(x_{nk})^* Q_{p'}^n [(E_n(y_{nk}))_{k}]_k ) \\
  &=& \summ_{nk} tr(x_{nk}^*Q_{p'}^n[ (E_n(y_{nk}))_{k}]_k) \\
  &=& tr( (x_{nk})^* Q_{p'}[ ( E_n(y_{nk}))]   ) \pl. \\[-1.675cm]
 \mel \qed

We want to remove the additional assumption that $M$ has
separable predual.\hhz

\begin{lemma}\label{sep} Let $M$ be a von Neumann algebra, $\phi$
a normal faithful state, $E_n$ a sequence of normal
conditional expectations onto von Neumann
subalgebras $N_n$. For every separable subalgebra
$A\subset M$ there exists a subalgebra $A\subset
B\subset M$ with separable dual, for all $\nen$
$E_n(B)\subset B$, and a normal conditional
expectation $\tilde{E}:M\to B$ such that $\phi \circ
\tilde{E}=\phi$.
\end{lemma}\hhz

{\bf Proof:} The proof is a modification of \cite[Corollary
3.5]{Kir}. Let $A_1$ be separable, $^*$-closed (but not
necessarily norm closed) algebra such that for all $x\in A$,
 \[ \sup_{y\in N, \noo y\rrm\le 1} |\phi(yx)| \lel
   \sup_{y\in A_1, \noo y\rrm\le 1}  |\phi(yx)|  \pl \]
and moreover $E_n(A)\subset A_1$ for all $n\in \nz$. Repeating
 this process, we obtain a separable $^*$-closed
subalgebra $A_\8=\bigcup_k A_k$ such that the embedding
$i:L_1(A_\8,\phi)\to L_1(N,\phi)$, $i(x.\phi)=x.\phi$ is
isometric. The dual map $\tilde{E}=i^*:M\to M$ is a normal
conditional expectation onto the $\si$-weak operator closure $B$
of $A_\8$. By construction, we have $\phi \circ \tilde{E}=\phi$
and
 $E_n(A_\8) \subset A_\8\cap N_n$. Since, $E_n$ is $\si$-strongly continuous,
$E_n(B)$ is a von Neumann subalgebra  of $N_n\cap B$ for all
$\nen$.\qed

\begin{samepage}
\begin{theorem}\label{dual} Let $1\le p \le \8$,
then $L_p(M,(E_n);\ell_2^C)$ is a Banach space and
the Cauchy-Schwarz inequality \ref{CS} $i)$ holds.
Let $1<p \, , \, p'<\8$ with
$\frac1p+\frac{1}{p'}=1$ and $\gamma_p$ the constant
from \ref{Stein}. If the sequence $(N_n)$ is either
increasing or decreasing  and $(x_n)$ a sequence in
$L_p(M,(E_n);\ell_2^C)$, then
 \for
 \lefteqn{ \noo (x_n)\rrm_{L_p(M,(E_n);\ell_2^C)}}
 \\
  & & \kl \gamma_p \pl \sup \left\{ \bigg | \summ_n tr(D^{\frac{1}{p'}}b_n^*x_n) \bigg |
  \pl \Bigg  |  \pl  \noo
  (b_nD^{\frac{1}{p'}} )\rrm_{L_{p'}(M,(E_n);\ell_2^C)} \kl 1\right\} \pl .
 \mel
\end{theorem}\hhz\end{samepage}

{\bf Proof:} Since the triangle inequality and  the
Cauchy-Schwarz inequality are checked for two
sequences, it suffices to consider a countably
generated subalgebra $B$ of $M$. Indeed,  according
to Lemma \ref{sep} we can even assume that
$E_n(B)\subset B$ and there exists a
$\phi$-preserving conditional expectation
$\tilde{E}:M\to B$. Then $L_{\frac{p}{2}}(B)$,  is
a  subspace of $L_{\frac{p}{2}}(M)$.  Hence,
$L_p(B,(E_n);\ell_2^C)$ is a subspace of
$L_p(M,(E_n);\ell_2^C)$ and hence there is no loss
of generality to assume that $M_*$ is separable.
Then  first  assertions follow from Proposition
\ref{CS}. For the last inequality, we   apply Lemma
\ref{dual0} to the image
$u_p(L_p(M,(E_n);\ell_2^C))$, where $u_p$ is the
isometric embedding from Proposition \ref{CS}.
Indeed, the projection $R_pST_p$ satisfies the
assumptions of Lemma \ref{dual0} and the norm is
bounded by the universal constant $\gamma_p$ from
Theorem \ref{Stein}. Finally, we note that $u_p$
preserves the duality bracket
 \for
 tr(u_{p'}[(b_nD^{\frac{1}{p'}})]^*u_p[(a_nD^{\frac1p})])
 &=&  \summ_n tr(D^{\frac{1}{p'}}E_n(b_n^*a_n) D^{\frac{1}{p}})
 \lel    \summ_n \phi( E_n(b_n^*a_n)) \\
 &=&   \summ_n \phi( b_n^*a_n) \lel
   \summ_n tr(D^{\frac{1}{p'}} b_n^*a_n D^{\frac1p}) \pl .
   \\
 \mel
Therefore the duality formula from Lemma \ref{dual0} is also
valid for $L_p(M,(E_n);\ell_2^C)$. \qed

\setcounter{equation}{0}

\section{\bf The dual version of Doob's inequality for \boldmath $1\le p\le 2$ \unboldmath}

In this section, we start with  an elementary proof of the dual
version of Doob's inequality for $p=2$ and show how the complex
interpolation method can be used to extend the inequality to  the
interval $1\le p\le 2$. Then we provide the duality argument which
justifies the name 'dual version of Doob's inequality'. In the
following, we consider a normal faithful state $\phi$, a von
Neumann algebra $N$  and a sequence of  von Neumann subalgebras
$N_n$ with conditional expectations $E_n:N\to N$ such that
 \[ \phi \circ E_n \lel \phi \]
for all $\nen$. (Note that in case of a tracial state such
conditional expectations always exist \cite{TAK}.) Let us stress
that in addition to the last section we also always assume that
the sequence $N_n$ is \underline{increasing}.\hhz

\begin{lemma}\label{p2} Let
  $(x_n)$ be a sequence of positive
elements in $L_2(N)$, then
 \[ \noo \summ_n E_n(x_n)\rrm_2 \kl 2 \noo \summ_n x_n \rrm_2 \pl
 .\]
\end{lemma}\hhz

{\bf Proof:} By monotonicity, $(1.3)$ and Corollary \ref{Ap3},  it
suffices to prove this inequality for finite sequences. Using
Proposition \ref{baby}, Lemma \ref{mod1} and positivity as in
$(1.3)$, see    \cite[Proposition 33]{Ter}, we deduce from
H\"older's inequality
 \for
  \noo \summ_n E_n(x_n)\rrm_2^2 &=& \summ_{nk}
 tr(E_n(x_n)E_k(x_k)) \\
 &=& \summ_{n\le k} tr(E_n(x_n)E_k(x_k)) \pl +\pl
    \summ_{n> k} tr(E_n(x_n)E_k(x_k)) \\
 &=& \summ_{n\le k} tr(E_k(E_n(x_n)x_k)) \pl +\pl
    \summ_{n> k} tr(E_n(x_n E_k(x_k)) ) \\
 &=& \summ_{n\le k} tr(E_n(x_n)x_k) \pl +\pl
    \summ_{n> k} tr(x_n E_k(x_k) ) \\
 &=& \summ_{k} tr ( \left ( \summ_{n\le k} E_n(x_n)\right ) x_k) \pl +\pl
    \summ_{n} tr(x_n \left ( \summ_{k<n} E_k(x_k)\right ) ) \\
 &\le& \summ_{k} tr ( \left (  \summ_{n} E_n(x_n)\right) x_k) \pl +\pl
    \summ_{n} tr(x_n \left ( \summ_{k} E_k(x_k)\right ) ) \\
 &=& 2 tr(\left(\summ_n x_n\right)  \left (\summ_k E_k(x_k)\right ) )
 \kl  2 \noo \summ_n x_n \rrm_2 \pl \noo \summ_n E_n(x_n)
 \rrm_2 \pl .
 \mel
Hence, we get

\begin{samepage}\for
 \noo \summ_n E_n(x_n)\rrm_2 &\le&   2 \noo \summ_n
 x_n\rrm_2 \pl. \\[-1.65cm]
 \mel \qed\end{samepage}

\hz


\begin{lemma}\label{intp1} If  $(DD_2)$ holds with constant $c_2$
and   $1\le p\le 2$, then  for  all sequences $(x_n)$ and $(y_n)$
in $L_{2p}(N)$
 \[ \noo \summ_n E_n(x_n^*y_n)\rrm_p  \kl c_2^{\frac{2(p-1)}{p}} \pl  \noo \summ_n x_n^*x_n
  \rrm_p^{\frac 12}
  \noo \summ_n y_n^*y_n \rrm_p^{\frac 12}   \pl
  .\]
\end{lemma}\hhz

{\bf Proof:} Let us first prove the assertion for finite sequences
and  $p=2$ or $p=1$.  We start with $(DD_1)$.  Using Lemma
\ref{mod1}, we get
  \for
  \noo \summ_n E_n(x_n^*x_n)\rrm_1 &=& tr( \summ_n
  E_n(x_n^*x_n)) \lel  \summ_n tr(E_n(x_n^*x_n))\\
  &=& \summ_n tr(x_n^*x_n) \lel    \noo \summ_n x_n^*x_n \rrm_1   \pl .
 \mel
By the density of elements of the form $x_n=a_nD^{\frac1p}$, the
Cauchy-Schwarz inequality, see Theorem \ref{dual}, implies with
Lemma \ref{p2} for $p\in \{1,2\}$
 \begin{eqnarray}
 \begin{minipage}{13cm} \vspace*{-0.5cm}
 \for
 \noo \summ_n E_n(x_n^*y_n)\rrm_p  &\le&
 \noo \summ_n E_n(x_n^*x_n)\rrm_p^{\frac  12 } \pl
   \noo \summ_n E_n(y_n^*y_n)\rrm_p^{\frac  12 }
   \\ &\le&  c_p \pl
    \noo \summ_n x_n^*x_n\rrm_p^{\frac  12 } \pl
 \noo \summ_n y_n^*y_n \rrm_p^{\frac  12 } \pl ,\mel\\[-0.9cm]
 \end{minipage}
 \end{eqnarray}
where $c_2$  is the constant given in the assumption and $c_1=1$.
To use interpolation, we consider finite sequences $(x_n)$ and
$(y_n)$ such that
\[ \noo \summ_n x_n^*x_n \rrm_p \lel 1 \lel  \noo \summ_n
 y_n^*y_n\rrm_p \pl .\]
We define  $X\lel \sum_n x_n^*x_n$, $Y\lel \sum_n y_n^*y_n$. Their
support projections  are denoted by $q_X$ and $q_Y$ and are in
$N$, see  \cite[Proposition 4. 2) c), Proposition 12]{Ter}. Since
$X^{-\frac 12}q_X$, $q_YY^{-\frac 12}$ are well-defined unbounded
operators, we can define
 \[ v_n \lel x_n  X^{-\frac 12}q_X
  \quad ,\quad w_n \lel y_nY^{-\frac 12} q_Y \]
Note that $x_n^*x_n\le X$ and  $y_n^*y_n\le Y$ implies
$x_n=x_nq_X$, $y_n=y_nq_Y$, respectively and according to Lemma
\ref{fact0} we have  $v_n\in N$, $w_n\in N$. Then, we observe
 \for
  \summ_n v_n^*v_n &=& q_XX^{-\frac 12}  X    X^{-\frac 12} q_X \le
  q_X\le 1 \pl ,\\
  \summ_n w_n^*w_n &=& q_Y Y^{-\frac 12} Y    Y^{-\frac 12} q_Y \le
  q_Y\le 1 \pl .
  \mel
Let $\theta$ be determined by 
$\frac1p=\frac{1-\theta}{1}+\frac{\theta}{2}$. According to 
Kosaki's interpolation  \cite{Kos} 
 \[  [L_{2,L}(N,\phi),L_{4,L}(N,\phi)]_{\theta}\lel L_{2p,L}(N,\phi)\quad \mbox{and} \quad 
  [L_{2,R}(N,\phi),L_{4,R}(N,\phi)]_{\theta}\lel L_{2p,R}(N,\phi) \] 
with respect to the state $\phi(x)=tr(Dx)$. By approximation, we 
may assume that  there are continuous  functions $X(z)$, $Y(z)$ 
on the strip $\{0\le Re(z)\le 1\}$ with values in $N$, analytic 
in the interior, such that 
$X^{\frac12}=D^{\frac{1}{2p}}X(\theta)$, 
$Y^{\frac12}=Y(\theta)D^{\frac{1}{2p}}$ and 
 \for 
  \sup_{t} \max\{\|D^{\frac{1}{2}} X(it)\|_2, \| D^{\frac14} X(1+it)\|_4\} &\le&   1 ,\\
  \sup_{t} \max\{\| Y(it)D^{\frac12}\|_{2}, \|  Y(1+it)D^{\frac14}\|_4\} &\le&   1 .
  \mel 
Now, we note Kosaki's  symmetric interpolation result 
\[ [L_{1,sym}(N,\phi),L_{2,sym}(N,\phi)]_\theta \lel L_{p,sym}(N,\phi)\pl. \] 
Hence, by $(3.1)$ and H\"older's inequality we get 
 \for 
 \lefteqn{  \noo \summ_n E_n(X^{\frac12}v_n^*w_nY^{\frac12})\rrm_p 
 \lel   \noo \summ_n D^{\frac{1}{2p}} 
 E_n(X(\theta)v_n^*w_nY(\theta))D^{\frac{1}{2p}}
  \rrm_p  }\\
  & &\kl  \sup_t   \noo \summ_n D^{\frac{1}{2}} E_n(X(it)v_n^*w_nY(it))D^{\frac{1}{2}}
  \rrm_1^{1-\theta} \\
  & &  \pll\pll  \sup_{t}   \noo \summ_n D^{\frac{1}{4}} E_n(X(1+it)v_n^*w_nY(1+it))D^{\frac{1}{4}}
  \rrm_2^{\theta} \\
  & & \kl    \sup_t   \noo \summ_n D^{\frac{1}{2}} X(it)v_n^*v_n X(it)^*D^{\frac{1}{2}}\rrm_1^{\frac{1-\theta}{2}}  
   \noo \summ_n D^{\frac{1}{2}} Y(it)^*w_n^*w_n Y(it)^*D^{\frac{1}{2}}\rrm_1^{\frac{1-\theta}{2}}  \\
  & & \pll \pll c_2^\theta \pl    \sup_t   \noo \summ_n D^{\frac{1}{4}} X(1+it)v_n^*v_n X(1+it)^*D^{\frac{1}{4}}\rrm_1^{\frac{\theta}{2}}  \\
  & & \pll \pll \pll   \pll \pll  \sup_t   \noo \summ_n D^{\frac{1}{4}} Y(1+it)^*w_n^*w_n Y(1+it)D^{\frac{1}{4}}\rrm_2^{\frac{\theta}{2}}  \\
  & & \kl \sup_t  \noo D^{\frac{1}{2}} X(it)X(it)^*D^{\frac12}\rrm_1^{\frac{1-\theta}{2}}  
  \sup_t  \noo D^{\frac{1}{2}} Y(it)^*Y(it)D^{\frac12}\rrm_1^{\frac{1-\theta}{2}}  \\
 & &   \pll \pll \pll  \sup_t  \noo D^{\frac{1}{4}} X(1+it)X(1+it)^*D^{\frac14}\rrm_2^{\frac{1-\theta}{2}}  
 \pl c_2^{\theta} \pl  \sup_t  \noo D^{\frac{1}{4}} Y(1+it)^*Y(1+it)D^{\frac14}\rrm_2^{\frac{1-\theta}{2}}  \\
 & & \kl  c_2^{\theta} \lel c_2^{\frac{2(p-1)}{p}} \pl . 
    \mel
The assertion is proved.\qed 

{\bf Proof of Theorem \ref{ep} in  the case $1\le p\le 2$:} For
$1\le p\le 2$ and a sequence of positive elements $(z_n)\subset
L_p(N)$, we can apply Lemma \ref{intp1} to $x_n=y_n=z_n^{\frac
12}$ and deduce the assertion for the sequence $(z_n)$.\qed

The  duality argument relies on the following norm  for  sequences
$(x_n)\subset L_p(N)$
 \[ \noo (x_n) \rrm_{\V_p(N;\ell_1)} \lel \inf  \noo
 \summ_{nj}
 v_{nj}v_{nj}^*\rrm_p^{\frac 12} \pl \noo \summ_{nj} w_{jn}^*w_{jn}\rrm_p^{\frac
 12} \pl .\]
Here the infimum is taken over all (double indexed) sequences
$(v_{nj})$ and $(w_{nj})$  such that for all $n$
 \[ x_n \lel \summ_j v_{nj}w_{jn} \pl .\]
We require norm convergence for $p<\8$ and convergence in the
$\si$-weak operator topology for $p=\8$. In fact, we think of
$x_n$ being obtained by matrix multiplication of a row with a
column vector. We denote by $\V_p(N;\ell_1)$ the set of all
sequences admitting a decomposition as above.\hhz

\begin{rem}{\rm This norm is motivated by the following
characterization of a normal,  decomposable  map $T:\ell_\8\to N$,
see \cite{Pau}.  Indeed, a normal map is decomposable if and only
if there are sequences $(x_n)\subset N$, $(y_n)\subset N$ such
that $T(e_n)=y_nx_n$ and}
\[ \noo \summ_n y_ny_n^*\rrm_N \pl \noo \summ_n x_n^*x_n\rrm_N \pl
<\pl \8 \pl .\]
\end{rem}\hhz

\begin{lemma} \label{CSA} If $(DD_p)$ holds, then
 \for
 \noo \summ_n E_n(x_n^*y_n)\rrm_p  \kl c_p \pl
    \noo \summ_n x_n^*x_n\rrm_p^{\frac  12 } \pl
 \noo \summ_n y_n^*y_n \rrm_p^{\frac  12 } \pl  ,
 \mel
The linear map $T:\V_p(N;\ell_1)\to \V_p(N;\ell_1)$,
$T((x_n))=(E_n(x_n))$  satisfies
 \[ \noo T\rrm \kl c_p  \pl .\]
\end{lemma}\hhz

{\bf Proof:} We have seen in $(3.1)$ that the first inequality is
 an immediate consequence of the Cauchy-Schwarz
inequality, see Theorem \ref{dual}. As for the second assertion,
we can assume that $N_*$ is separable and use the Kasparov maps
$u_p^n:L_p(N,E_n) \to L_p(N_n;\ell_2^C)$ from Proposition
\ref{CS1}. Let $x_n=\sum_j v_{nj}w_{nj}$, then
 \for
  \noo (E_n(x_n)) \rrm_{\V_p(N;\ell_1)} &=& \noo (\summ_{nj} u_p^n(v_{nj}^*)^*u_p^n(w_{nj})) \rrm_{\V_p(N;\ell_1)} \\
  &\le& \noo \summ_{nj}  u_p^n(v_{nj}^*)^*u_p^n(v_{nj}^*) \rrm_p^{\frac12} \pl
   \noo \summ_{nj}  u_p^n(w_{nj})^*u_p^n(w_{nj}) \rrm_p^{\frac12} \pl \\
  &=& \noo \summ_{nj} E_n(v_{nj}v_{nj}^*) \rrm_p^{\frac12} \pl
   \noo \summ_{nj} E_n(w_{nj}^*w_{nj}) \rrm_p^{\frac12}  \\
  &\le& c_p \pl  \noo \summ_{nj} v_{nj}v_{nj}^* \rrm_p^{\frac12}
 \pl     \noo \summ_{nj} w_{nj}^*w_{nj} \rrm_p^{\frac12}  \pl .
 \mel
Taking the infimum over all these decompositions, we obtain the
assertion.\qed

Let us state some elementary properties of the space
$\V_p(N;\ell_1)$.\hhz

\begin{lemma} \label{elepop} For $1\le p\le \8$ the set $\V_p(N;\ell_1)$ is a Banach
space. For $1\le p <\8$, the set  $\V_p^0$ of elements
\[ x_n \lel \summ_j v_{nj}w_{jn} \]
such that $card\{(j,n) | v_{nj}\neq 0 \mbox{ or } w_{jn}\neq
0\}<\8$ is dense in $\V_p(N;\ell_1)$. If $(x_n)$ is a sequence
 then
  \[ \noo \summ_n x_n \rrm_p \kl \noo (x_n) \rrm_{\V_p(N;\ell_1)} \]
and equality holds if all the $x_n$'s are positive.
\end{lemma}\hhz

{\bf Proof:} The proof of the triangle inequality is completely
elementary, see also \cite{Pvp},  but essential. Indeed, let
$\eps>0$ and
 \for
  x_n &=&  \summ_{j_1} v_{nj_1}w_{j_1n} \quad \mbox{and} \quad
  y_n \lel \summ_{j_2} v_{nj_2}w_{j_2n} \pl ,
  \mel
such that
 \for
  \noo \summ_{nj_1} v_{nj_1}v_{nj_1}^*\rrm_p &=&  \noo \summ_{nj_1}
  w_{j_1n}^*w_{j_1n}\rrm_p \kl (1+\eps) \noo
  (x_n)\rrm_{\V_p(N;\ell_1)} \pl ,\\
  \noo \summ_{nj_2} v_{nj_2}v_{nj_2}^*\rrm_p &=&  \noo \summ_{nj_2}
  w_{j_2n}^*w_{j_2n}\rrm_p \kl (1+\eps) \pl \noo
  (y_n)\rrm_{\V_p(N;\ell_1)} \pl .
 \mel
We have
 \[ x_n+y_n \lel \summ_{j_1} v_{nj_1}w_{j_1n} +  \summ_{j_2}
 v_{nj_2}w_{j_2n}\]
and the triangle inequality in $L_p(N)$ implies
 \for
 \noo \summ_{nj_1} v_{nj_1}v_{nj_1}^* +  \summ_{nj_2}
 v_{nj_2}v_{nj_2}^* \rrm_p
 &\le& \noo \summ_{nj_1} v_{nj_1}v_{nj_1}^*\rrm_p + \noo
 \summ_{nj_2}
 v_{nj_2}v_{nj_2}^* \rrm_p\\
 &\le& (1+\eps) (\noo
  (x_n) \rrm_{\V_p(N;\ell_1)} +\noo
  (y_n) \rrm_{\V_p(N;\ell_1)}) \pl .
 \mel
Similarly,
 \[ \noo \summ_{nj_1} w_{j_1n}^*v_{j_1n} +  \summ_{nj_2}
 w_{j_2n}^*w_{j_2n} \rrm_p \kl (1+\eps) (\noo
   (x_n) \rrm_{\V_p(N;\ell_1)} +\noo
  (y_n) \rrm_{\V_p(N;\ell_1)}) \pl \]
and the assertion follows with $\eps \to 0$.  We consider the
spaces of column matrices, row matrices, $L_p(N;\ell_2^C(\nz^2))$,
$L_p(N;\ell_2^R(\nz^2)) \subset L_p(B(\ell_2(\nz^2))\bar{\ten}
N)$, respectively. Using,
 \[ \Phi( (x_{nk})\ten  (y_{nk})) \lel (\summ_k x_{nk}y_{nk})_{\nen}
 \pl ,
 \]
we deduce that $\V_p(N;\ell_1)$ is isomorphic to a quotient space
of the projective  tensor product $L_p(N;\ell_2^R(\nz^2))
\ten_{\pi}
 L_p(N;\ell_2^C(\nz^2))$,
see \cite{DF} for a definition and basic properties of the
projective tensor product. Hence  $\V_p(N;\ell_1)$ is complete.
The image under $\Phi$ of pairs of finite sequences generates
$\V_p^0$. According to Corollary \ref{Ap3} finite sequences are
dense in the column and row spaces and hence $\V_p^0$ is dense in
$\V_p(N;\ell_1)$. Finally, let $(x_n)$ be a sequence of positive
elements. Clearly, $x_n=x_n^{\frac12}x_n^{\frac12}$ and hence
 \[ \noo (x_n) \rrm_{\V_p(N;\ell_1)} \kl \noo \summ_n x_n \rrm_p
 \pl .\]
On the other hand if $x_n=\sum_j v_{nj}w_{nj}$, we deduce from
H\"older's inequality
 \for
  \noo \summ_n x_n \rrm_p &=& \noo \summ_{nj} (e_{1,nj} \ten v_{nj}) (e_{nj,1}
  \ten w_{nj}) \rrm_p \\
  &\le&   \noo \summ_{nj} e_{1,nj}\ten  v_{nj}\rrm_{2p}
  \noo \summ_{nj} e_{nj,1}\ten  w_{nj}\rrm_{2p} \\
  &=&  \noo \summ_{nj} v_{nj}v_{nj}^*\rrm_p^{\frac12}
   \noo \summ_{nj} w_{nj}^*w_{nj}\rrm_p^{\frac12} \pl.
  \mel
Taking the infimum yields the assertion.\qed

Inspired by Pisier's vector-valued $L_p(N,\tau;\ell_\8)$ space,
we define for $0<p\le \8$
 \[ \noo \sup_n |x_n| \rrm_p \lel  \noo (x_n)\rrm_{\L_p(N;\ell_\8)} \lel  \inf_{x_n=ay_nb}
 \noo a\rrm_{2p} \pl  \noo b\rrm_{2p} \pl \sup_n \noo y_n\rrm_N  \pl ,\]
where the infimum is taken over all $a,b\in L_{2p}(N)$ and all
bounded sequences $(y_k)$. If $N$ is a hyperfinite, finite von
Neumann algebra, this space coincides with $L_p(N,\tau;\ell_\8)$
in the sense of Pisier \cite{Pvp}. The first (formal) notation is
suggestive and facilitates understanding our inequalities  in
view of the commutative theory.  For positive elements, we will
drop the absolute value.  Let us note that Haagerup's work
\cite{Haai} shows that the equality $\L_1(N;\ell_\8)=
L_1(N)\ten_{\wedge} \ell_\8$ (operator space projective tensor
product) only holds for injective von Neumann algebras. However,
this does not affect the following factorization result which is,
nowadays, a standard application of the Grothendieck-Pietsch
version  of the Hahn-Banach theorem, see \cite{Plf,Pvp}.\hhz

\begin{prop}\label{N} Let $1\le p<\8$. If $p=1$, then
$\V_1(N;\ell_1)=\ell_1(L_1(N))$ holds with  equal norms. If $1< p
\, , \, p'<\8$ satisfy $\frac 1p+\frac{1}{p'}=1$, then
 \[  \V_p(N;\ell_1)^*\lel \V_{p'}(N;\ell_\8) \]
holds isometrically.
\end{prop}\hhz

{\bf Proof:} If $z_n=ay_nb$ and $(y_n)$ is a  bounded sequence, we
deduce from H\"older's inequality  for all $(v_{nj})\subset
L_{2p}(N)$, $(w_{jn})\subset L_{2p}(N)$
 \for
 \lefteqn{\bet \summ_n tr(z_n\summ_j v_{nj} w_{jn} ) \rag \lel
  \bet  \summ_{nj} tr(ay_nbv_{nj}w_{jn})\rag
  \lel
   \bet \summ_{nj} tr(y_nbv_{nj}w_{jn}a)\rag  }\\
 & & \kl  \sup_n \noo y_n\rrm_\8  \pl \summ_{nj} \noo
  bv_{nj}w_{jn}a \rrm_1 \\
  & & \kl  \sup_n \noo y_n\rrm_\8 \pl  \kla \summ_{nj} \noo bv_{nj}\rrm_2^2 \mer^{\frac12} \pl
         \kla \summ_{nj} \noo w_{jn}a\rrm_2^2 \mer^{\frac12} \\
  & & \lel  \sup_n \noo y_n\rrm_\8 \pl  \noo b\summ_{nj}v_{nj}v_{jn}^*b^*\rrm_1^{\frac 12} \pl
      \noo a^*\summ_{nj} w_{nj}^*w_{jn}a\rrm_1^{\frac12} \\
  & & \kl  \sup_n \noo y_n\rrm_\8 \pl \noo b\rrm_{2p'} \pl  \noo \summ_{nj}v_{nj}v_{nj}^*\rrm_p^{\frac 12} \pl
          \noo a \rrm_{2{p'}} \pl       \noo \summ_{nj} w_{jn}^*w_{jn}\rrm_p^{\frac12}
          \pl .
 \mel
Hence, $\V_{p'}(N;\ell_\8)\subset \V_p(N;\ell_1)^*$. Using
$\ell_1(L_1(N))\subset \V_1(N;\ell_1)$, we deduce the equality
$\ell_1(L_1(N))= \V_1(N;\ell_1)$. Now we show that  for $1< p
<\8$ all the functionals are in $\V_{p'}(N;\ell_\8)$. Let
$\psi:\V_p(N;\ell_1)\to \cz$ be a norm one functional.
 Using
$\ell_1(L_p(N))\subset \V_p(N;\ell_1)$, we can assume that there
exists a sequence  $(z_n) \subset L_{p'}(N)$ such that
 \for
 \psi[(x_n)] \lel \psi_{(z_n)}[ (x_n)] \lel \summ_n tr(z_nx_n) \pl .
  \mel
Let us denote by $B=B_{L_{p'}(N)}^+$ the positive part of the
unit ball in $L_{p'}(N)$. $B$ is compact when  equipped with the
$\si(L_{p'}(N),L_{p}(N))$-topology. The definition of
$\V_p(N;\ell_1)$ implies with the geometric/arithmetric mean
inequality
 \for \bet \summ_{nj}
 tr(z_nv_{nj}w_{jn}) \rag &=&
 \bet \psi[ (\summ_j v_{nj}w_{jn})_n ] \rag \kl
    \noo \summ_{nj} v_{nj}v_{nj}^*\rrm_p^{\frac 12 } \pl
       \noo \summ_{nj} w_{jn}^*w_{jn}\rrm_p^{\frac 12}\\
 &\le & \frac{1}{2}\pl \sup_{c,d\in B}
 [ \summ_{nj} tr(v_{nj}v_{nj}^*c)+
  \summ_{nj} tr(w_{jn}^*w_{jn}d)] \pl .
 \mel
Since the right hand side remains unchanged under multiplication
with signs $\eps_{nj}$, we deduce
 \for
 \summ_{nj}
|tr(z_nv_{nj}w_{jn})| &\le& \frac{1}{2} \pl
  \pl \sup_{c,d\in B}
 [ \summ_{nj} tr(v_{nj}v_{nj}^*c)+
  \summ_{nj} tr(w_{jn}^*w_{jn}d)] \pl .
  \mel
Following  the  Grothendieck-Pietsch separation
argument as in \cite{Plf}, we observe that the $C$
given by  the functions
\[   f_{v,w}(c,d) \lel \summ_n [tr(v_nv_n^*c)+
 tr(w_n^*w_nd)-2|tr(z_nv_nw_n)|] \]
is disjoint from the cone $C_{-}=\{g | \sup g
<0\}$. Here $v=(v_n)$ and $w=(w_n)$ are finite
sequences and hence $f_{v,w}$ is continuous with
respect to the product topology on $B\times B$.
Since $f_{v,w}+f_{\tilde{v},\tilde{w}}$ can be
obtained by taking the $(\tilde{v}_n,\tilde{w}_n)$'s
to the right of the finite sequence $(v_n,w_n)$, we
deduce that $C$ is a cone. Hence, there exist a
measure $\mu$ on $B\times B$ and a scalar $t$ such
that for all $g\in C_-$ and $f\in C$
\[ \intt_{B\times B} g \pl d\mu \pl <t\kl \intt_{B\times B} f \pl d\mu \pl .\]
Since we are dealing with cones, it turns out that $t=0$ and
$\mu$ is positive.  Therefore,   we can and will assume that
$\mu$ is  a probability measure.  We define  the positive
elements $c$ and $d$ by their projections
\[ a \lel \intt_{B\times B}  c \pl  \pl d\mu(c,d)\pll ,
 \pll
   b  \lel \intt_{B\times B} d  \pl  \pl d\mu(c,d)\pl . \]
By convexity of $B$, we deduce $a,b\in B$. Hence, we obtain
 \for
 \summ_n 2 \pl |tr(z_nv_nw_n)| &\le& \intt_{B\times B} \summ_n
 [tr(v_nv_n^*c) +  tr(w_nv_n^*d)] \pl d\mu(c,d) \\
 &=& \summ_n \intt_{B\times B} tr(v_nv_n^*c) \pl d\mu(c,d) +
              \intt_{B\times B} tr(w_n^*w_nd) \pl d\mu(c,d)\\
 &=& \summ_n [tr(v_nv_n^*a) + tr(w_n^*w_nb)] \pl .
 \mel
Using once more $2st \lel \inf_{r>0} (rs)^2+(r^{-1}t)^2$, we get
 \begin{eqnarray}
 \begin{minipage}{13cm} \vspace*{-0.5cm}
 \for
 \summ_n |tr(z_nv_nw_n)|  &\le& \kla \summ_n tr(v_nv_n^*a)
 \mer^{\frac 12} \pl \kla \summ_n tr(w_n^*w_nb) \mer^{\frac 12}
 \\
 &=& \kla \summ_n \noo a^{\frac 12}v_n\rrm_2^2 \mer^{\frac 12} \pl
     \kla \summ_n \noo b^{\frac12}w_n^*\rrm_2^2 \mer^{\frac12} \pl.
 \mel \end{minipage}
 \end{eqnarray}
Let $q_a$, $q_b \in N$  be the support projections
of $a$, $b$, respectively.  Consider,
$d_a=a^{p'}+(1-q_a)D(1-q_a)$, $D$ the density of
$\phi$. Then $\phi_{d_a}(x)=tr(d_ax)$ is a normal,
faithful state on $N$ and according to Lemma
\ref{Ap0} $d_a^{\frac12}N$ is dense in $L_2(N)$.
Hence,
 \[ q_a d_a^{\frac12}N \lel a^{\frac{p'}{2}} N \lel
 a^{\frac12}a^{\frac{p'}{2p}}N \subset a^{\frac 12}L_{2p}(N) \]
shows that $a^{\frac 12}L_{2p}(N)$ is dense in $q_aL_2(N)$.
Similarly, $b^{\frac 12}L_{2p}(N)$ is dense  in $q_bL_2(N)$ and
therefore $(3.2)$ implies that for every $\nen$ there is a
contraction $T_n:q_aL_2(N)\to q_bL_2(N)$ such that for all $v,w\in
L_{2p}(N)$
 \[ tr(wz_nv) \lel  (b^{\frac12}w^*,T_n(a^{\frac12}v)) \lel
 tr(wb^{\frac12}T_n(a^{\frac12}v)) \pl .\]
This means $T_n$ is a bounded extension of the
densely defined  hermitian form
  \[ (b^{-\frac12}q_b(h'),z_na^{-\frac12}q_a(h)) \lel
  (b^{-\frac12}q_bh',z_na^{-\frac12}q_ah)  \]
Using the density of $a^{\frac12}L_{2p}(N)$ and
$b^{\frac12}L_{2p}(N)$ it is easily checked that
$q_bT_nq_a$ is affiliated with $N$. Since $T_n$ is
bounded, we deduce $q_bT_nq_a\in N$. On the other
hand, we have for $v\in L_{2p}(N)$ and $w\in
L_{2p}(N)$
\[ |tr(z_n(1-q_a)vw)| \kl tr(a(1-q_a)vv^*)^{\frac12}
 tr(w^*wb)^{\frac12} \lel 0 \pl \]
and
\[ |tr(z_nvw(1-q_b)| \kl tr(avv^*)^{\frac12}
 tr(w^*w(1-q_b)b)^{\frac12} \lel 0 \pl .\]
This shows  $z_n\lel q_bz_nq_a$ and therefore
 \[ z_n \lel q_bz_nq_a \lel q_bb^{\frac12} q_bb^{-\frac 12} z_n a^{-\frac 12}q_a a^{\frac
 12}q_a \lel b^{\frac12}y_n a^{\frac12} \pl .\]
The assertion is proved because  $\V_p^0$ is dense in
$\V_p(N;\ell_1)$ and hence the  functional $\psi$ is uniquely
determined by the sequence $(z_n)$.\qed

\begin{rem}\label{Np} Let  $1\le p<\8$ and
$(z_n)\subset L_{p'}(N)$ a sequence of  positive elements, then
 \[ \noo \sup_n z_n \rrm_{p'} \lel \sup \left\{ \bet \summ_n tr(z_nx_n)
 \rag \mitt x_n\ge 0\pl ,\pl \noo \summ_n x_n\rrm_p \le 1\right\}
 \pl .\]
Moreover, there exists a positive element $a\in L_{2p'}(N)$ and a
sequence of positive elements $y_n$ such that
 \[ z_n \lel ay_na \quad \mbox{and} \quad \noo a\rrm_{2p'}^2  \pl \sup_n \noo y_n\rrm_\8 \lel
  \noo \sup_n z_n \rrm_{p'} \pl .\]
For positive elements $(x_n) \subset L_p(N)_+$, we also have
 \[ \noo \summ_n x_n \rrm_p \lel \sup \left\{ \bet \summ_n
 tr(x_nz_n)\rag \mitt z_n \ge 0\pl ,\pl  \noo (z_n)\rrm_{\V_{p'}(N;\ell_\8)} \kl 1\right\}  \]
and therefore the cones  of positive sequences in
$\V_p(N;\ell_1)$ respectively $\V_{p'}(N;\ell_\8)$  are in
duality.
\end{rem}

{\bf Proof:}  For positive elements $(z_n)$ satisfying
 \for
  \bet \summ_n tr(z_nx_n)\rag &\le& \noo \summ_n x_n \rrm_p \pl ,
 \mel
for all sequences   of positive elements $(x_n)\subset L_p(N)_+$,
we deduce
 \for
 \summ_{j,n} |tr(z_nv_{nj}v_{nj}^*)| &\le&
 \bet \summ_{n} tr(z_n( \summ_j v_{nj}v_{nj}^*))\rag
 \kl   \noo \summ_{nj} v_{nj}v_{nj}^*\rrm_p \\
 &=&  \sup_{c\in B}
     \summ_{j,n} tr(v_{nj}v_{nj}^*c) \pl .
 \mel
Using  the Hahn-Banach separation argument in the space of
continuous functions on $B$, we  obtain a positive element $a$ in
the unit ball of $L_{p'}(N)$ such that
\[ \summ_n |tr(z_nv_nv_n^*)| \kl \summ_n tr(v_nv_n^*a) \pl .\]
Since $a^{-\frac12}z_na^{-\frac12}$ is  positive, this inequality
still ensures that all the $y_n=a^{-\frac12}z_na^{-\frac12}$ are
positive contractions in $N$. The last equality follows
immediately from Lemma \ref{elepop} and the duality between the
positive parts of $L_p(N)$ and $L_{p'}(N)$. \qed

\begin{rem} Proposition \ref{N} and Remark \ref{Np} can easily be
modified for uncountable (ordered) index sets by requiring the
inequality for all countable (ordered) subsets or for an essential
supremum. This is  helpful in the context of continuous
filtrations.
\end{rem}

The required duality argument is now very simple.\hhz

\begin{lemma}\label{eck}  Let $1< p \le  \8$ and  $\frac1p+\frac{1}{p'}=1$
assume that $(DD_{p'})$ holds with constant $c_{p'}$, then for
every $y\in L_p(N)$
 \[ \noo \sup_n |E_n(y)|   \rrm_{p} \kl c_{p'}    \pl \noo y\rrm_p\pl .\]
Moreover, for every sequence of positive elements $(y_n)$
 \[ \noo \sup_n |\summ_{j\ge n} E_n(y_j)|  \rrm_{p} \kl c_{p'} \pl \noo
 \summ_n y_n \rrm_p \pl .\]
\end{lemma}\hhz

{\bf Proof:} Indeed, as observed in Lemma  \ref{CSA} and using
Lemma \ref{elepop}, we deduce that  $(DD_{p'})$ implies that the
linear map  $T:\V_{p'}(N;\ell_1)\to L_p(N)$, $T((x_n))=\sum_n
E_n(x_n)$ satisfies $\noo T\rrm\le c_{p'}$. By duality  and
Proposition \ref{N}, we deduce for all $y\in L_p(N)$
  \[ \noo \sup_n |E_n(y)| \rrm_{p} \lel \noo
  T^*(y)\rrm_{\V_{p'}(N;\ell_1)^*} \kl \noo T^*\rrm \noo y\rrm_p
  \kl c_{p'}\noo y\rrm_p \pl .\]
Given a sequence   of positive elements $(y_n)$, we consider
$y=\sum_j y_j$ and obtain
 \[ \noo \sup_n |E_n(y)|  \rrm_{p}  \kl c_{p'} \pl \noo y\rrm_p\lel c_{p'} \pl \noo \summ_j
 y_j\rrm_p \pl .\]
However, for positive elements $(x_n)\subset L_{p'}(N)$, we deduce
by positivity
 \for
  \summ_n tr(E_n(\summ_{j\ge n} y_j)x_n) &\le & \summ_n  \summ_j
  tr(E_n(y_j)x_n) \kl \noo (E_n(y))_{\nen}  \rrm_{\V_{p}(N;\ell_\8)} \noo
  \summ_n x_n\rrm_{p'} \pl .
  \mel
Hence, Remark \ref{Np} implies
 \for
  \noo \sup_n \summ_{j\ge n} E_n(y_j)\rrm_{p}
  &\le&
 \noo \sup_n E_n(y) )_{\nen}  \rrm_{p}  \kl c_{p'} \pl \noo y\rrm_p \lel c_{p'} \pl\noo \summ_j
 y_j\rrm_p \pl .
  \mel
The assertion is proved.\qed

\begin{theorem}\label{doob}
For  $1<p\le \8$ there exists a constant $c_p$ such that for every
sequence $(N_n)$  of  von Neumann subalgebras with sequence of
$\phi$-invariant conditional expectations $(E_n)$ satisfying
$E_{n}E_m=E_{\min(n,m)}$  and for  every $x\in L_p(N)$ there
exist $a,b\in L_{2p}(N)$ and a bounded sequence $(y_n)\subset N$
such that
 \[ E_n(x) \lel ay_nb  \quad \mbox{and} \quad \noo a \rrm_{2p} \noo b\rrm_{2p} \sup_n \noo y_n\rrm_\8 \kl c_p \noo x\rrm_p  \pl .\]
If  $x$ is positive, one can in addition assume  that $b=a^*$ and
all the $y_n$'s  are positive.
\end{theorem}\hhz

{\bf Proof for  $2\le p\le \8$:} This follows immediately from
Lemma \ref{eck}, Lemma \ref{intp1} and Lemma \ref{p2}. Using that
$E_n(x)$ is positive for positive $x$, the addition follows from
Remark \ref{Np}. \qed

\section{\bf The dual version of Doob's inequality for \boldmath$2\le p <\8$\unboldmath}

In  our approach to $(DD_p)$ in the range $2\le p<\8$, our aim is
to obtain the same kind of inequalities for the maximal function
as in Garsia's book \cite{Gar}. As mentioned in the introduction,
we are forced to use more duality arguments because $0\le a\le b$
implies  $a^\beta\le b^\beta$ only for $0\le \beta\le 1$ and
therefore most of the elementary proofs in Garsia's  book  are no
longer valid in the non-commutative case. We will make the same
assumptions about $N$, $(N_n)$, $(E_n)$ and $\phi$, $D$ as in the
previous section. In particular, the sequence $N_n$ is supposed
to be  increasing. \hhz

\begin{lemma}\label{al}
Let $0\le x\le z \in L_1(N)$ such that the support projection of
$z$ is $1$. Let  $1\le \al \le 2$, then
 \[ tr(z^{\frac{1-\al}{2}}(z^\al-x^\al)z^{\frac{1-\al}{2}}) \kl 2 tr(z-x) \pl .\]
\end{lemma}\hhz

{\bf Proof:} We define $\beta=\al-1\in(0,1)$ and observe that
$0\le x\le z$ implies
 \[  x^{1-\beta} \kl z^{1-\beta}  \pl .\]
We apply Lemma \ref{fact0} to $x^{\frac \beta 2 }$ and $z^{\frac
\beta 2}$ and deduce from $x^{\beta}\le z^{\beta}$ that $v=
x^{\frac{\beta}{2}} z^{-\frac{\beta}{2}}$ is a contraction in $N$.
In particular $v^*=z^{-\frac{\beta}{2}}x^{\frac{\beta}{2}}\in N$
and
 \[ a \lel z^{-\frac{\beta}{2}}x^{\frac{1+\beta}{2}}  \lel v^*x^{\frac12} \in L_2(N) \pl .\]
For all elements $a,b\in L_2(N)$, we note that
 $(a-b)^*(a-b)\ge 0$ implies with tracial property of the trace
\[ tr(a^*b)+tr(b^*a)  \kl tr(a^*a) + tr(b^*b) \lel tr(a^*a) + tr(bb^*)  \pl .\]
We define $b=
z^{\frac{\beta}{2}}x^{\frac12-\frac{\beta}{2}}$
which is in $L_2(N)$ by H\"older's inequality.
Then, we observe
 \for
   tr(a^*b) &=&  tr(x^{\frac{1+\beta}{2}} z^{-\frac{\beta}{2}} z^{\frac{\beta}{2}} x^{\frac12-\frac{\beta}{2}} ) \lel tr(x)\pl , \\
   tr(b^*a) &=&   tr(x^{\frac12-\frac{\beta}{2}}   z^{\frac{\beta}{2}} z^{-\frac{\beta}{2}} x^{\frac{1+\beta}{2}}) \lel tr(x)\pl .
 \mel
Hence, we deduce
 \for
  2tr(x) &=& tr(a^*b)+tr(b^*a) \kl    tr(a^*a) +tr(bb^*) \\
  &=& tr(x^{\frac{1+\beta}{2}}  z^{-\frac{\beta}{2}} z^{-\frac{\beta}{2}}x^{\frac{1+\beta}{2}})
   +  tr(z^{\frac{\beta}{2}}x^{\frac12-\frac{\beta}{2}}x^{\frac12-\frac{\beta}{2}}z^{\frac{\beta}{2}}) \\
  &=& tr(x^{\frac{1+\beta}{2}}z^{-\beta} x^{\frac{1+\beta}{2}})     + tr(z^{\frac{\beta}{2}}x^{1-\beta}z^{\frac{\beta}{2}}) \\
  &\le&  tr(x^{\frac{1+\beta}{2}}z^{-\beta} x^{\frac{1+\beta}{2}})  + tr(z^{\frac{\beta}{2}}z^{1-\beta}z^{\frac{\beta}{2}}) \\
  &=&  tr(x^{\frac{\al}{2}}z^{1-\al}x^{\frac{\al}{2}} ) + tr(z) \pl .
   \mel
This implies
 \[ -tr(x^{\frac{\al}{2}}z^{1-\al}x^{\frac{\al}{2}} )  +tr(z)   \kl -2tr(x) +tr(z)+tr(z)  \lel 2 tr(z-x)  \pl .\]
The assertion follows from $z^{\frac{1-\al}{2}}x^\al
z^{\frac{1-\al}{2}}=aa^*\in L_1(N)$ and
 \for
 tr(x^{\frac{\al}{2}}z^{1-\al}x^{\frac{\al}{2}} )
 &=&  tr(a^*a)  \lel  tr(aa^*) \lel   tr(z^{\frac{1-\al}{2}}x^\al z^{\frac{1-\al}{2}}) \pl . \\[-1.65cm]
 \mel \qed

The following proposition is a modification of the
Theorem in the appendix of Pisier, Xu's paper
\cite{PX} and enables us to apply duality.\hhz

\begin{prop}\label{dual2}   Let $1\le r'< 2< r\le \8$, and
$\frac{ 1}{r'}+\frac1r \lel 1$, then for all $(x_j)\subset
L_r(N)$ and $(y_n)\subset L_{r'}(N,(E_n);\ell_2^C)$
 \[ \bet \summ_n tr(y_n^*x_n) \rag
 \kl \sqrt{2} \pl \noo (y_n) \rrm_{L_{r'}(N,(E_n);\ell_2^C)}
 \pl
 \noo \summ_j x_j^*x_j \rrm_{\frac{r}{2}}^{\frac 12}  \pl .\]
Moreover,
 \[ \bet \summ_n tr(y_n^*x_n) \rag
 \kl \sqrt{2} \pl \noo (y_n) \rrm_{L_{r'}(N,(E_n);\ell_2^C)}
 \pl
 \noo \sup_n \summ_{j\ge n} E_n(x_j^*x_j) \rrm_{\frac{r}{2}}^{\frac 12}  \pl .\]
\end{prop}\hhz

{\bf Proof:} Let us assume that both sequences are finite, i.e.
$x_j=0=y_j$ for $j\ge m$.  By density,   we can moreover assume
that $y_j=a_jD^{\frac{1}{r'}}$ with $a_j\in N$ and
 \[  \noo  \summ_{j=1}^{m} D^{\frac{1}{r'}}E_j(a_j^*a_j)D^{\frac{1}{r'}} \rrm_{\frac{r'}{2}} \pl <\pl  1 \pl .\]
By continuity, we can assume  that there is an $\eps>0$ such that
 \[  \noo \eps D^{\frac{2}{r'}}+ \summ_{j=1}^{m} D^{\frac{1}{r'}}E_j(a_j^*a_j)D^{\frac{1}{r'}} \rrm_{\frac{r'}{2}} \kl 1 \pl .\]
Let $1\le q\le \8$ such that $\frac{1}{q}+\frac{2}{r}=1$. For
$\nen$, we define
 \[ S_n \lel \kla  \eps D^{\frac{2}{r'}} +\summ_{j=1}^n D^{\frac{1}{r'}}E_j(a_j^*a_j)D^{\frac{1}{r'}} \mer^{\frac{r'}{2q}}  \in  L_q(N_n) \subset L_q(N) \pl .\]
The support projection of $S_n$ is $1$  and $\eps
D^{\frac1q}\le  S_n$. According to Lemma
\ref{fact0}, we deduce that
$w_n:=D^{\frac{1}{2q}}S_n^{-\frac12}\in N_n$. Hence
 \[ y_nS_n^{-\frac12} \lel a_nD^{\frac{1}{r'}}S_n^{-\frac12} \lel
  a_nD^{\frac12} D^{\frac{1}{2q}}S_n^{-\frac12}\lel a_n
  D^{\frac12}w_n  \in L_2(N) \pl . \]
In particular, $S_n^{-\frac12}y_n^*\in L_2(N)$ and
 \[  S_n^{-\frac12}y_n^*y_nS_n^{-\frac12}
  \lel  w_n^*D^{\frac12}a_n^*a_nD^{\frac12}w_n  \in L_1(N) \pl .\]
Moreover,
 \for
 E_n(S_n^{-\frac12}y_n^*y_nS_n^{-\frac12}) &=&
 E_n(w_n^*D^{\frac12}a_n^*a_nD^{\frac12}w_n)
 \lel   w_n^*D^{\frac12}E_n(a_n^*a_n)D^{\frac12}w_n  \\
 &=& S_n^{-\frac12} D^{\frac{1}{r'}} E_n(a_n^*a_n)
 D^{\frac{1}{r'}}S_n^{-\frac12} \pl .
 \mel
This implies with  the Cauchy-Schwarz inequality
 \for
 \bet \summ_n tr(y_n^*x_n) \rag &=& \bet \summ_n tr(x_ny_n^*) \rag
 \lel
 \bet \summ_n tr((x_nS_n^{\frac12})(S_n^{-\frac12} y_n^*))
 \rag\\
 &=& \bet \summ_n tr(S_n^{-\frac12} y_n^*x_nS_n^{\frac12} ) \rag \\
 &\le& \kla \summ_n tr(S_n^{-\frac12} y_n^*y_nS_n^{-\frac12}) \mer^{\frac12}  \pl
       \kla \summ_n tr( x_n^*x_nS_n) \mer^{\frac12}  \pl   \\
 &=& \kla \summ_n tr(E_n(S_n^{-\frac12} y_n^*y_nS_n^{-\frac12})) \mer^{\frac12}  \pl
       \kla \summ_n tr( x_n^*x_nS_n) \mer^{\frac12}  \pl   \\
 &=& \kla \summ_n tr(S_n^{-\frac12}D^{\frac{1}{r'}} E_n(a_n^*a_n)D^{\frac{1}{r'}}S_n^{-\frac12}) \mer^{\frac12}  \pl
       \kla \summ_n tr( E_n(x_n^*x_n)S_n) \mer^{\frac12}  \pl  .
  \mel
To estimate the first term, we define $\al=\frac{2}{r'}\in [1,2]$
and notice that
 \[  1-\al \lel 1-\frac{2}{r'} \lel 1-2+\frac{2}{r} \lel -\frac 1q \pl .\]
For fixed $n$, we define $x=S_{n-1}^q$ and  $z=S_n^q$. Since
$\frac{r'}{2}\le 1$, we have
 \[ x \lel \kla \eps D^{\frac{2}{r'}}+ \summ_{j=1}^{n-1} D^{\frac{1}{r'}}E_j(a_j^*a_j)D^{\frac{1}{r'}} \mer^{\frac{r'}{2}}
 \kl  \kla \eps D^{\frac{2}{r'}}+ \summ_{j=1}^{n} D^{\frac{1}{r'}}E_j(a_j^*a_j)D^{\frac{1}{r'}} \mer^{\frac{r'}{2}} \lel z  \pl .\]
Then, we note  that $z^{\frac{1-\al}{2}} \lel z^{-\frac{1}{2q}}
\lel S_n^{-\frac12}$. Hence Lemma \ref{al} implies
 \for
  tr(S_n^{-\frac12} D^{\frac{1}{r'}} E_n(a_n^*a_n)D^{\frac{1}{r'}}S_n^{-\frac 12}) &=&
  tr(z^{\frac{1-\al}{2}} (z^{\al}-x^{\al})z^{\frac{1-\al}{2}} )
  \kl  2 tr(z-x) \\
  &=&  2 tr(S_n^q-S_{n-1}^q) \pl .
   \mel
Therefore, we obtain
 \for
 \summ_n tr(S_n^{-\frac12} D^{\frac{1}{r'}} E_n(a_n^*a_n)D^{\frac{1}{r'}}S_n^{-\frac12} )
 &=&   \summ_n 2 tr( S_n^q-S_{n-1}^q) \lel  2 tr(S_m^q) \\
 &=&    2 \noo \eps D^{\frac{2}{r'}}  + \summ_{j=1}^m  D^{\frac{2}{r'}} E_j(a_j^*a_j)D^{\frac{2}{r'}}\rrm_{\frac{r'}{2}}^{\frac{r'}{2}} \kl 2 \pl .
 \mel
Now, we want to estimate the second term. Let us define
 \[ \theta_j \lel S_j-S_{j-1} \pl .\]
and note that $\theta_j\in L_q(N_j)$. As usual, we set
$S_{-1}=0$. Moreover, $r'\le 2q$ implies that $\theta_j$ is
positive. Then, we deduce with $E_jE_n=E_{min(j,n)}$
 \for
  \summ_n tr(E_n(x_n^*x_n)S_n) &=& \summ_{j\le n }
  tr(E_n(x_n^*x_n)\theta_j) \\
  &=& \summ_j tr( \summ_{n\ge j} E_n(x_n^*x_n)\theta_j) \\
  &=& \summ_j tr(E_j\bigg( \summ_{n\ge j} E_n(x_n^*x_n) \bigg) \theta_j) \\
  &=& \summ_j tr( \summ_{n\ge j} E_j(x_n^*x_n)\theta_j) \\
 \mel

Now, we can continue in two different ways
 \for
    \summ_j tr( \summ_{n\ge j} E_j(x_n^*x_n)\theta_j)
  &\le&     \noo \sup_j E_j(\summ_{n\ge j}  x_n^*x_n) \rrm_{\frac r2}  \pl   \noo \summ_j \theta_j \rrm_q \\
  &=&    \noo \sup_j E_j(\summ_{n\ge j}  x_n^*x_n) \rrm_{\frac r2} \pl
    \noo S_m\rrm_q \\
  &\le&   \noo \sup_j E_j(\summ_{n\ge j}  x_n^*x_n) \rrm_{\frac r2} \pl .
  \mel
By homogeneity, we obtain the second assertion
 \for
  \bet \summ_n tr(y_n^*x_n) \rag &\le&
  \sqrt{2} \noo    \sup_j E_j(\summ_{n\ge j}  x_n^*x_n) \rrm_{\frac r2}^{\frac12}
  \noo (y_n) \rrm_{L_{r'}(N,(E_n);\ell_2^C)} \pl .
 \mel
Using for all $j\in \nz$ that $E_j(\theta_j)=\theta_j$, we also
get by positivity
 \for   \summ_j tr( \summ_{n\ge j} E_j(x_n^*x_n)\theta_j)
 &=& \summ_j tr( \summ_{n\ge j} x_n^*x_n E_j(\theta_j)) \\
 &=& \summ_j tr( \summ_{n\ge j} x_n^*x_n \theta_j ) \\
 &\le& tr( \summ_{n} x_n^*x_n \summ_j \theta_j) \\
 &\le& \noo \summ_n x_n^*x_n\rrm_{\frac{r}{2}} \pl \noo \summ_j\theta_j \rrm_q \\
 &\le& \noo \summ_n x_n^*x_n\rrm_{\frac{r}{2}} \pl \noo S_m\rrm_q \kl
 \noo \summ_n x_n^*x_n\rrm_{\frac{r}{2}}  \pl .
 \mel
Again by homogeneity, we deduce
 \for
  \bet \summ_n tr(y_n^*x_n) \rag &\le&
  \sqrt{2} \noo  \summ_{n}  x_n^*x_n  \rrm_{\frac r2}^{\frac12}
  \noo (y_n) \rrm_{L_{r'}(N,(E_n);\ell_2^C)} \pl .\\[-1.675cm]
 \mel\qed

\begin{rem} For $r'=r=2$ the assertion is trivially true because
$L_{2}(N,(E_n);\ell_2^C)=\ell_2(L_2(N))$ and
 \for
 \summ_n \noo x_n\rrm_2^2 &=& \noo \summ_n x_n^*x_n \rrm_1 \lel
 tr(\summ_n x_n^*x_n) \lel  tr(E_1(\summ_n x_n^*x_n))\\
 &=& \noo E_1(\summ_n x_n^*x_n) \rrm_1 \kl \noo \sup_j
 E_j(\summ_{n\ge j} x_n^*x_n)\rrm_1 \pl .
 \mel
Proposition \ref{dual2}  also shows that the $BMO_C$ and $H_1^C$
duality from \cite{PX} is still valid in the non tracial  case.
\end{rem}\hhz

\begin{cor}\label{kick} Let $1\le q <\8$ and
$\gamma_{2q}$  the constant in Stein's inequality, then for all
$(x_n)\subset L_{2q}(N)$
  \[ \noo \summ_n E_n(x_n^*x_n)\rrm_{q} \kl 2 \gamma_{2q}^2 \pl
  \noo \sup_n E_n(\summ_{j\ge n} x_j^*x_j)\rrm_{q} \pl .\]
\end{cor}\hhz

{\bf Proof:} We define $r=2q$ and note
 \[ \noo \summ_n E_n(x_n^*x_n)\rrm_{r}^{\frac12}  \lel
    \noo (x_n)\rrm_{L_{p}(N,(E_n);\ell_2^C)}  \pl .\]
Therefore the assertion follows from Proposition \ref{dual2} and
Theorem \ref{dual}.\qed

{\bf Proof of \boldmath $(DD_p)$ \unboldmath  for $1< p<\8$:} We
define $r=2p>2$  Let $(z_n)\subset L_{p}(N)$ be a sequence of
positive elements and define $x_n=z_n^{\frac12}$. Hence by
Theorem \ref{dual} and Proposition \ref{dual2}, we deduce
 \for
 \lefteqn{  \noo \summ_n E_n(z_n)\rrm_p^{\frac12}  \lel   \noo \summ_n
  E_n(x_n^*x_n)\rrm_p^{\frac12}}\\
  & & \kl   \gamma_{r}  \pl \sup_{\noo (y_n) \rrm_{L_{r'}(N,(E_n);\ell_2^C)}\le 1}    \bet \summ_n tr(y_n^*x_n) \rag \\
  & & \kl  \gamma_{r} \pl \sqrt{2}\pl  \noo \summ_n x_n^*x_n\rrm_p^{\frac12} \pl .
  \mel
The assertion follows with $c_p\le 2 \gamma_{2p}^2$. \qed

{\bf Proof of Theorem \ref{0.2.} and \ref{doob} for $1<p\le 2$:}
This follows from $(DD_p)$ via Lemma \ref{eck} and Remark
\ref{Np}.\qed

\begin{rem} Let $N_*$ be separable and $\psi$ be a functional on
$L_{p'}(N,(E_n);\ell_2^C)$, then there exists a sequence $(x_n)$
such that
 \[ \psi( (y_n)) \lel \summ_n tr(x_n^*y_n) \quad \mbox{and} \quad
 \noo \sup_n E_n(\summ_{j\ge n} x_j^*x_j)
 \rrm_{\frac{p}{2}}^{\frac12} \kl d_{\frac{p}{2}} \noo \psi\rrm_{L_{p'}(N,(E_n);\ell_2^C)^*}  \pl .\]
Here $d_{\frac{p}{2}}$ is the constant in Doob's inequality from
Theorem \ref{0.2.}. The assertion  yields an extension of the
$BMO_C$-$H_1^C$ duality for $2< p<\8$ and  fails for $p=2$. The
proof uses the Kasparov isomorphism from  Proposition \ref{CS}.
We leave it to the interested reader.
\end{rem}

Answering a question by G. Pisier, we can even produce an
asymmetric version of Doob's inequality. Indeed,  let $1<p\le \8$
and consider $1\le r,s < \8$ such that
$\frac{2}{p'}=\frac{1}{s}+\frac{1}{t}$.  Given $x\in L_p(N)$ and
sequences $(v_{nj})\subset L_{r'}(N)$, $(w_{jn})\subset
L_{s'}(N)$, we deduce   from the Cauchy-Schwartz inequality
\ref{Stein}, $(DD_{s})$   and  $(DD_{t})$ that
 \for
  \bet \summ_{n,j} tr(E_n(x)v_{nj}w_{jn})  \rag &=&
  \bet \summ_{n,j} tr(xE_n(v_{nj}w_{jn})) \rag \\
  &\le& \noo x\rrm_p \pl \noo \summ_{nj} E_n(v_{nj}w_{jn}) \rrm_{p'} \\
  &\le& \noo x\rrm_p \pl \noo \summ_{nj} E_n(v_{nj}v_{nj}^*) \rrm_{s}^{\frac12} \pl
  \noo \summ_{nj} E_n(w_{jn}^*w_{jn}) \rrm_{t}^{\frac12}  \\
  &\le&  c_{s}^{\frac12} c_{t}^{\frac12} \pl
  \noo x\rrm_p \pl \noo \summ_{nj} v_{nj}v_{nj}^*  \rrm_{s}^{\frac12} \pl
  \noo \summ_{nj} w_{jn}^*w_{jn} \rrm_{t}^{\frac12}   \pl .
  \mel
Similar as in Proposition \ref{N}, we deduce the existence of
bounded a sequence $(z_n)$ and elements $a\in L_{2s'}(N)$,
$b\in L_{2t'}(N)$ such that  $E_n(x)=bz_na$ and
 \[ \noo a \rrm_{2 s'} \noo b\rrm_{2 t'}  \sup_n \noo z_n\rrm_\8 \pl  \kl
       c_{s}^{\frac12} c_{t}^{\frac12} \pl   \noo x\rrm_p \pl .\]
Note that $\frac{1}{2t'}+\frac{1}{2s'}=\frac1p$ and therefore we
have proved the following asymmetric version of Theorem
\ref{0.2.}. (The assumption $2<q,r$ is indeed necessary
\cite{DJ}.)

\begin{cor} Let $1<p\le \8$ and $2<q,r \le \8$ such that $\frac{1}{q}+\frac{1}{r}=\frac1p$. Then for every $x\in L_p(N)$ there
exists a sequence $(z_n)\subset N$ and $a\in L_{q}(N)$, $b\in
L_{r}(N)$ such that $E_n(x)=az_nb$ and
 \[ \noo a\rrm_q \noo b\rrm_r \sup_n \noo z_n\rrm_\8    \kl c(p,q,r) \pl \noo x\rrm_p \pl .\]
\end{cor}

\section{\bf Applications}

In this section, we present first applications of Doob's
inequality in terms of submartingales, Doob decomposition. We
make the same assumptions about $N$, $(N_n)$, $(E_n)$, $\phi$ and
$D$ as in the previous section and start with almost immediate
consequences of the dual version of Doob's inequality.\hhz

\begin{cor} \label{sub1}
Let $1<  p\le\8$ and $(z_n)$ be an adapted sequence of positive
elements, i.e. $z_n\in L_p(N_n)_{+}$. If for all $\nen$
 \[ z_n \le E_n(z_{n+1})  \quad \mbox{and}  \quad
 \sup_m \noo z_m\rrm_p \pl <\pl \8\pl ,\]
then there exist a positive element $a\in L_{2p}(N)_+$ and    a
sequence of positive contractions $(y_n)\subset N$ such that
 \[ z_n \lel  ay_na \quad \mbox{and}  \quad
 \noo a\rrm_{2p}^2 \kl  c_{p'} \pl \sup_m \noo z_m\rrm_p \pl
 .\]
\end{cor}\hhz

{\bf Proof:} It suffices to consider $p<\8$. We note that for
$n\le m$
 \[ z_n\le E_n(z_{n+1}) \kl E_n(E_{n+1}(z_{n+2})) \lel E_n(z_{n+2})
 \kl \cdots \kl  E_n(z_m) \pl .\]
Let $\frac1p+\frac{1}{p'}=1$. In order to estimate the norm in
$\L_{p}(N;\ell_\8)$, we refer to Remark \ref{Np}. Given a finite
sequence $(x_n)$ of positive elements such that $x_n=0$ for $n\ge
m$, we deduce from $(DD_{p'})$
 \for
 \bet \summ_n tr(z_nx_n) \rag &=& \summ_n tr(E_n(z_nx_n)) \lel
  \summ_n tr(E_n(z_n)E_n(x_n)) \\
 &\le &   \summ_n tr(E_n(z_m)E_n(x_n))  \lel
 \summ_n tr(E_n(z_mE_n(x_n)))\\
  &=&    \summ_n tr(z_m E_n(x_n))  \kl       tr(z_m\summ_n E_n(x_n))  \\
  &\le&  \noo z_m\rrm_p \pl \noo \summ_n E_n(x_n)\rrm_{p'} \kl
  c_{p'} \pl \noo z_m\rrm_p \pl \noo \summ_n x_n \rrm_{p'}
 \mel
Hence, the assertion follows from Remark \ref{Np}.\qed

\begin{cor} Let $2< p\le \8$ and $(z_n)$ be an adapted seqeunce,
i.e. $z_n\in L_p(N_n)$. If for all $n\in \nz$
\[ z_{n}^*z_{n}\kl E_{n}(z_{n+1}^*z^{}_{n+1}) \pl ,\]
then there exist a positive element $a\in L_{p}(N)$
and a sequence $(y_n)$ of contractions  such that
 \[  z_n\lel y_n a \quad \mbox{and} \quad
 \noo a\rrm_p \kl c_{\frac{p}{p-2}}^{\frac12}    \pl \sup_m \noo z_m\rrm_p \pl .\]
\end{cor}\hhz

{\bf Proof:} We apply Corollary \ref{sub1} to $z_n'=z_n^*z_n\in
L_{\frac{p}{2}}(N_n)$ and obtain positive contractions
$(v_n)\subset N$ and $a\in L_p(N)$ such that
 \[ z_n^*z_n\lel av_na \quad \mbox{and} \quad   \noo a\rrm_p^2 \kl c_{\frac{p}{p-2}}  \pl \sup_m \noo z_m^*z_m\rrm_{\frac p2} \pl .\]
If $q_a$ is the support projection of $a$, we see that $y_n \lel
v_n^{\frac{1}{2}}a^{-1}q_a$ satisfies the assertion. \qed

Let us mention an  immediate application of $(DD_p)$ in terms  of
the Doob decomposition of the square function.  Given a
martingale sequence $x=\sum_k d_k(x)$, where
$d_k(x)=E_k(x)-E_{k-1}(x)$ (and $E_0(x)=0$), we recall that one
of the square functions is given by
  \[ s_c(x) \lel \summ_k d_k(x)^*d_k(x) \]
The  square function is the discrete analogue of the quadratic
variation term see \cite{BS,Ko} for more details. The Doob
decomposition of $s_c(x)$ is given by  the martingale part
 \[ V_n(x) \lel \summ_{k=1}^n d_k(x)^*d_k(x)-E_{k-1}(d_k(x)^*d_k(x)) \]
and  the predictable part
 \[ W_n(x) \lel \summ_{k=2}^n E_{k-1}(d_k(x)^*d_k(x)) \in L_p(\M_{k-1}) \pl . \]
Note that
 \[ V_n(x) + W_n(x) \lel \summ_{k=1}^n d_k(x)^*d_k(x) \pl .\]

\begin{cor} Let $2<p<\8$ and $x \in L_p(N)$, then
 \[     \sup_n \max\{ \noo W_n(x)\rrm_{\frac p2}^{\frac12}, \noo  V_n(x) \rrm_{\frac{p}{2}}^{\frac12}\}
   \kl  \al_p (1+c_{\frac p2})^{\frac12}  \noo x\rrm_p  \pl .\]
Here $\al_p$ is an absolute constant. In particular, there exists
an element $V_\8(x)$ in $L_{\frac p2}(N)$ such that
$V_n(x)=E_n(V_\8(x))$.
\end{cor}\hhz

{\bf Proof:} Using $(DD_{\frac p2})$ for the sequence
$(E_{k-1})$  and the non-commutative Burkholder-Gundy inequality
\cite{PX,JX}, we deduce
 \for
 \noo W_n(x) \rrm_{\frac{p}{2}} &\le& c_{\frac{p}{2}} \noo \summ_{k=2}^n d_k(x)^*d_k(x) \rrm_{\frac p2}
 \kl  c_{\frac{p}{2}}  \al_p^2 \noo x\rrm_p^2 \pl .
 \mel
Hence the triangle implies
 \for
  \noo V_n(x)\rrm_{\frac p2} &\le& \noo W_n(x)\rrm_{\frac{p}{2}} +\noo \summ_{k=1}^n d_k(x)^*d_k(x)\rrm_{\frac p2}   \kl
 (c_{\frac{p}{2}}  \al_p^2+1) \pl \noo x\rrm_p^2  \pl .
 \mel
By uniform convexity of $L_{\frac{p}{2}}(N)$, we obtain the limit
value  $V_\8(x)=\lim_n V_n(x)$ with the desired properties.
\qed

The next application yields norm estimates for $\sum_n p_n
E_n(x)q_n$ with respect to a sequence $(p_n)$, $(q_n)$ of
disjoint projections. This corresponds to  a double sided
non-adapted stopping time.\hhz

\begin{cor}
Let $1<p\le \8$,  $(v_n)$, $(w_n)$ be sequences of bounded
elements, then for all $x\in L_p(N)$
 \for
  \noo \summ_n v_n E_n(x)w_n \rrm_p &\le&
 c_{p'}  \pl \noo x \rrm_p  \pl
 \max\left\{ \noo \summ_n v_nv_n^* \rrm_\8^{\frac12},
 \noo \summ_n v_n^*v_n \rrm_\8^{\frac12} \right\} \\
 & & \pll \pl \pll \pl  \pll \pl   \pll \pl    \max\left\{ \noo \summ_n w_nw_n^* \rrm_\8^{\frac12},
 \noo \summ_n w_n^*w_n \rrm_\8^{\frac12} \right\}  \pl .
 \mel
\end{cor}\hhz

{\bf Proof:} The case $p=\8$ is obvious. Hence, we assume
$1<p<\8$. Let $y\in L_{p'}(N)$ and choose $y_1\in L_{2p'}(N)$,
$y_2\in L_{2p'}(N)$ such that $y=y_1y_2$ and
\[ \noo y\rrm_{p'} \lel \noo y_1\rrm_{2p'}^2 \lel \noo y_2\rrm_{2p'}^2  \pl .\]
Then, according to Lemma \ref{CSA}, we deduce from $(DD_{p'})$:
  \for
 \bet  tr(\summ_n v_n E_n(x)w_ny) \rag
  &=& \bet \summ_n tr(xE_n(w_nyv_n)) \rag
  \kl  \noo x\rrm_p \pl \noo \summ_n E_n(w_ny_1y_2v_n) \rrm_{p'} \\
  &\le& \noo x\rrm_p c_{p'} \pl
  \noo \summ_n w_ny_1y_1^*w_n^*\rrm_{p'}^{\frac 12} \pl
  \noo \summ_n v_n^*y_2^*y_2w_n\rrm_{p'}^{\frac 12} \pl .
 \mel
To conclude, we use \cite[Lemma 1.1]{PX}, see also \cite{JX},
 \for
 \noo \summ_n w_ny_1y_1^*w_n^*\rrm_{p'} &\le& \noo
 y_1y_1^*\rrm_{p'} \pl
 \max\left\{ \noo \summ_n w_nw_n^* \rrm_\8,
 \noo \summ_n w_n^*w_n \rrm_\8\right\} \\
 &=&   \noo y\rrm_{p'}
 \pl
 \max\left\{ \noo \summ_n w_nw_n^* \rrm_\8,
 \noo \summ_n w_n^*w_n \rrm_\8\right\} \pl.
 \mel
A similar argument applies for the other term and therefore taking
the supremum over all $y$ of norm 1 implies the assertion.\qed

\begin{rem} \label{scc} Let $1\le p<\8$. The non-commutative Doob inequality from Theorem
\ref{0.2.} implies the vector-valued Doob inequality in
$L_p(\Om,\Si,\mu;L_p(N))$.
\end{rem}\hhz

{\bf Proof:} It suffices to treat the discrete case. Let
$\N=L_\8(\Om,\Si,\mu)\bar{\ten}N$ and $(\Si_n)_{\nen}$ be an
increasing sequence of $\si$-subalgebras with conditional
expectations $(E_n)$. Let  $\N_n=L_\8(\Om,\Si,\mu)\bar{\ten}N$.
Then the conditional expectation $\E_n$ onto $\N_n$ is given by
$\E_n \lel E_n\ten id$. Let $f\in
L_p(\Om,\Si,\mu;L_p(N))=L_p(\N)$. According to Theorem \ref{0.2.}
there exist $a\in L_{2p}(\N), b\in L_{2p}(\N)$ and contractions
$(z_n)\subset \N$ such that
 \[ E_n \ten id_{L_p(N)}(f) \lel a z_n b  \quad \mbox{and}\quad \noo a\rrm_{2p} \noo b\rrm_{2p} \kl c_{p'} \noo f\rrm_p \pl .\]
Hence, for every $\om\in \Om$ and $\nen$
 \for
   \noo E_n(f)(\om) \rrm_{L_p(N)} &=& \noo a(\om)  z_n(\om)  b(\om)
  \rrm_{L_p(N)}
  \kl   \noo a(\om)\rrm_{L_{2p}(N)} \noo z(\om)\rrm_N
 \noo   b(\om)\rrm_{L_{2p}(N)} \\
 &\le&  \noo a(\om)\rrm_{L_{2p}(N)} \noo b(\om)\rrm_{L_{2p}(N)} \pl .
  \mel
H\"older's   inequality implies the assertion
 \for
 \kla  \intt_\Om \sup_n \noo E_n(f)(\om) \rrm_{L_p(N)}^p
 d\mu(\om)
 \mer^{\frac1p} &\le& \kla \intt_\Om  \noo a(\om)\rrm_{L_{2p}(N)}^p  \noo
 b(\om)\rrm_{L_{2p}(N)}^p d\mu(\om)  \mer^{\frac1p}      \\
 &\le& \noo a\rrm_{L_{2p}(\Om,\Si,\mu;L_{2p}(N)} \pl
 \noo b\rrm_{L_{2p}(\Om,\Si,\mu;L_{2p}(N)}\\
 &=& \noo a\rrm_{2p} \pl \noo b\rrm_{2p} \kl c_{p'} \noo f\rrm_p
 \pl .\\[-1.65cm]
 \mel\qed

In  the next application we want to relate  group actions with
$(DD_p)$. To illustrate this, we consider a finite von Neumann
algebra $N$ and an increasing sequence $(A_n)\subset N$ of finite
dimensional subalgebras with $1_N\in A_n$. Let $N_n=A_n'$ be the
relative commutant of $A_n$ in $N$. If $G_n$ denotes the unitary
group of $A_n$, we have a natural action $\al:G_n\to B(L_p(N))$
 \[ \al_n(u)(x) \lel uxu^* \]
such that the conditional expectation on the commutant $N_n$ is
given by
 \[ E_n(x) \lel E_{N_n}(x) \lel \intt_{G_n} uxu^* d\mu_n(u) \pl .\]
Let $G=\prod_n G_n$ and $\mu$ the product measure, then
  \[ \noo \summ_n E_n(x_n) \rrm_p \kl \kla \intt_G \noo \summ_n
  \al_n(u_n)(x_n) \rrm_p^p d\mu(u_1,u_2,..)  \mer^{\frac1p} \pl .\]
We will show that for a sequence $(x_n)$ of positive elements
even the right hand  side  can  be estimated by $\noo \sum_n
x_n\rrm_p$. For simplicity let us  use the random variables
$\al_n(x):G\to L_p(N)$, given for $\om=(u_1,u_2...)$ by
 \[ \al_n(x)(\om) \lel u_nxu_n^*  \pl .\]
For the special case of tensor products of finite dimensional von
Neumann algebras the following theorem implies $(DD_p)$.\hhz

\begin{theorem} Let $1<p<\8$ and $N$, $(A_n)$, $(N_n)$ be as above
and $(x_n)$ be a  sequence of positive elements, then
 \[ \kla \intt_G \noo \summ_n \al_n(x_n) \rrm_p^p d\mu\mer^{\frac1p}  \kl \kappa_p \pl  \noo \summ_n
 x_n \rrm_p \pl .\]
Here $\kappa_p$ is a constant which only depends on $p$.
\end{theorem}\hhz

{\bf Proof:} The assertion is obvious for $p=1$ and by
interpolation as  in Lemma \ref{intp1} it suffices to prove the
assertion for $p\ge 4$.  Let $(x_n)$ be a finite sequence of
positive elements.  Then, we observe
 \for
 \kla \intt_G   \noo \summ_n \al_n(x_n)\rrm_p^p d\mu   \mer^{\frac1p}
  &\le& \noo \summ_n E_n(x_n)\rrm_p + \kla \intt_G   \noo \summ_n \al_n(x_n)-E_n(x_n) \rrm_p^p  d\mu  \mer^{\frac1p}  \\
  &\le& c_{p} \noo \summ_n x_n \rrm_p + \kla \intt_G \noo \summ_n \al_n(x_n)-E_n(x_n) \rrm_p^p d\mu \mer^{\frac1p}  \pl .
 \mel
Let $\Si_n$ be the $\si$-algebra generated by the first $n$
coordinates in $G=\prod_k G_k$ and
$\N_n=L_\8(G,\Si_n,\mu)\bar{\ten} N\subset
L_\8(G,\Si,\mu)\bar{\ten} N$ with the corresponding conditional
expectation $\E_n$. Then, we note
 \[ \E_{n-1}(\al_n(x_n)) \lel \intt_{\prod_{k\ge n} G_k} u_n x_n u_n^*
 \pl d\mu_n(u_n)  d\mu_{n+1}(u_{n+1}) \cdots \lel E_n(x_n) \pl .\]
Hence the $n$-th martingale difference of $\sum_n \al_n(x_n)$
satisfies
 \[ d_n \lel \E_n(\sum_k \al_k(x_k))-\E_{n-1}(\sum_k \al_k(x_k)) \lel
  \al_n(x_n)-E_n(x_n)  \pl. \]
We apply the non-commutative Rosenthal inequality, see \cite{JX},
in this case and obtain
 \for
 \lefteqn{ \frac{1}{r_p}\pl  \kla \ez \noo \summ_n \al_n(x_n)-E_n(x_n) \rrm_p^p
 \mer^{\frac1p}}\\
 &\le& \kla \summ_n \noo  \al_n(x_n)-E_n(x_n)\rrm_p^p \mer^{\frac1p}
  + \noo \summ_n \E_{n-1}(d_nd_n^*+d_n^*d_n)\rrm_{\frac{p}{2}}^{\frac12} \\
 &\le& 2 \kla \summ_n \noo  x_n \rrm_p^p \mer^{\frac1p}
  + 2^{\frac12}  \noo \summ_n E_n(x_nx_n^*)+ E_n(x_n^*x_n) \rrm_{\frac{p}{2}}^{\frac12} \\
 &\le& 2 \kla \summ_n \noo  x_n \rrm_p^p \mer^{\frac1p}
  + 2^{\frac12}  c_{\frac p2}^{\frac12}  \noo \summ_n x_nx_n^*+x_n^*x_n  \rrm_{\frac{p}{2}}^{\frac12}
 \mel
Let $(\eps_n)$ be a sequence of independent Rademacher variables.
Using the triangle inequality and the orthogonality of the
$(\eps_n)$'s, we deduce as in \cite{LP}
 \for
 \max\left\{ \noo \summ_n x_n^*x_n  \rrm_{\frac{p}{2}}^{\frac12},
  \noo \summ_n x_nx_n^*  \rrm_{\frac{p}{2}}^{\frac12} \right \}
 &\le&  \kla \ez \noo \summ_n \eps_n x_n\rrm_p^2 \mer^{\frac12}
 \pl.
 \mel
By interpolation, we have
 \for \kla \summ_n \noo x_n \rrm_p^p
 \mer^{\frac1p}
 &\le& \noo \summ_n x_n^*x_n  \rrm_{\frac{p}{2}}^{\frac12}
 \pl .
 \mel
However, since the $x_n$ are positive, we deduce for any choice of
signs $\eps_n$ with positivity
 \for
  \noo \summ_n \eps_n x_n\rrm_p &\le& \noo \summ_{\eps_n=1} x_n\rrm_p +
  \noo \summ_{\eps_n=-1} x_n\rrm_p \kl 2 \pl \noo \summ_n x_n\rrm_p \pl .
 \mel
Hence, we obtain
 \for
 \kla \ez \noo \summ_n \al_n(x_n)-E_n(x_n) \rrm_p^p
 \mer^{\frac1p} &\le& r_p (2+ 4c_{\frac{p}{2}}^{\frac12})
 \noo \summ_n x_n \rrm_p \pl .
 \mel
The assertion is proved.\qed

\begin{rem}{\rm These methods can also be
used to show that  for every $f\in L_p(G;L_p(N))$ there exist
$a,b\in L_{2p}(\Om;L_{2p}(N))$ and a sequence of contractions
$(y_n)\subset L_\8(\Om)\bar{\ten} N$ such that
 \[ F_n(f) \lel a y_nb  \quad \mbox{and} \quad \noo a\rrm_{2p} \noo b\rrm_p \kl c_{p'}^2 \noo f\rrm_p \]
where
 \[ F_n(f)(g_1,g_2,...) \lel \a_n(g_n)f(g_1,...,g_n,..) \pl .\]
Note that the crossed product $N\rtimes_{(\al_n)} \prod G_n$ acts
on $L_p(\Om;L_p(N))$ and $F_n$ somehow removes the action of $G_n$
on $f$. Similar results hold for a von Neumann algebra with a
faithful normal state $\phi$ and  $\phi$-invariant, strongly
continuous group actions $\al_n:G_n\to Aut(N)$ of  compact groups
such that the centralizer algebras are increasing or decreasing.}
\end{rem}

\renewcommand{\baselinestretch}{1.0}
\newcommand{\sss}{\vspace{0cm}}

2000 Mathematics Subject Classification: 46L53, 46L52, 47L25.
\newline  Key-words: Non-commutative $L_p$-spaces, Doob's inequality

\end{document}